\documentclass{article}
\usepackage{amssymb,amsmath}
\usepackage{graphicx,enumerate}
\usepackage[utf8]{inputenc}
\usepackage[all]{xy}

\title{On pseudo-irreducibility and Boolean lifting property of filters in residuated lattices}

\author{Esmaeil Rostami \\ \footnotesize  Department of Pure Mathematics, Faculty
	of Mathematics and Computer,\\ \footnotesize Shahid
	Bahonar University of Kerman, Kerman,
	Iran\\ \footnotesize Email: e\underline{ }rostami@uk.ac.ir}

\date{}

\begin{document}
\maketitle

\begin{abstract}
In this paper, we introduce the notion of a pseudo-irreducible filter in a residuated lattice and compare this concept with related notions such as prime and maximal filters. Then, we recall the Boolean lifting property for filters and present useful characterizations for this property using pseudo-irreducible filters and the residuated lattice of fractions. Next, we study the Boolean lifting property of the radical of a filter. Furthermore, we introduce weak MTL-algebras and residuated lattices that have the transitional property of radicals decomposition (TPRD) as generalizations of several algebraic structures, including Boolean algebra, MV-algebra, BL-algebra, MTL-algebra, and Stonean residuated lattice. Moreover, by comparing weak MTL-algebras with other classes of residuated lattices, we address an open question concerning the Boolean lifting property of the radical of a residuated lattice. Finally, we give a topological answer to an open question about the Boolean lifting property of the radical of a residuated lattice. Several additional results are also obtained, further enriching the understanding of the Boolean lifting property in residuated lattices.

\end{abstract}

\textbf{Keywords}: Pseudo-irreducible filter, Boolean lifting property, Residuated lattice, Weak MTL-algebra, Transitional property of radicals decomposition (TPRD).
\newtheorem{definitie}{Definition}[section]
\newtheorem{propozitie}[definitie]{Proposition}
\newtheorem{remarca}[definitie]{Remark}
\newtheorem{exemplu}[definitie]{Example}
\newtheorem{intrebare}[definitie]{Open question}
\newtheorem{lema}[definitie]{Lemma}
\newtheorem{teorema}[definitie]{Theorem}
\newtheorem{corolar}[definitie]{Corollary}

\newenvironment{proof}{\noindent\textbf{Proof.}}{\hfill\rule{2mm}{2mm}\vspace*{5mm}}

\section*{Introduction}
Residuated lattices, introduced by Ward and Dilworth [25] in 1939 as a generalization of the lattice of ideals of rings, provide a unifying framework for studying algebraic structures arising in both algebra and logic. They include important classes of algebras, such as MV-algebras, MTL-algebras, and BL-algebras, which are central to fuzzy logic and non-classical  propositional logical, including $\L$ukasiewicz logic, H$\acute{a}$jek’s Basic Logic, and substructural logics, see [7] and [11].

In ring theory, the concept of comaximal ideal factorization is an important topic when considering rings in terms of ideal factorization. In this regard, pseudo-irreducible ideals are of particular importance, see [13] and [18]. As an extension, in [23], the notion of a pseudo-irreducible ideal in a De Morgan residuated lattice was introduced. Then the authors explored the relationships between such ideals and other important concepts, such as prime ideals and maximal ideals, in residuated lattices.   

On the other hand, one of the most useful techniques in ring theory is to first examine the properties of certain quotients and then \lq\lq lift\rq\rq  these properties to the original ring (see [22]). For example, \lq\lq lifting idempotents\rq\rq  is an instance of this technique. Building on this idea, the authors of [6] considered a lifting property for Boolean elements modulo the radical in MV-algebras and used this property to characterize maximal MV-algebras. Subsequently, this concept was generalized to BL-algebras and residuated lattices in [5], [9], [19], and [20]. The authors of [10] extended their study to the Boolean Lifting Property (BLP) modulo the radical in residuated lattices from an algebraic perspective. In [8], this concept was further extended to universal algebras, where they investigated lifting properties for idempotent and Boolean elements modulo filters. They also provided various algebraic and topological characterizations for the Boolean lifting property.

The Boolean lifting property raises numerous questions, making it a valuable subject for the effective study of residuated lattices. Moreover, its close connection to pseudo-irreducibility and the lifting of idempotents in commutative rings provides strong motivation for investigating these concepts in the context of residuated lattices. In this paper, we define the notion of a pseudo-irreducible filter in a residuated lattice and examine its relationship with the Boolean lifting property. We establish several new results in this direction, including characterizations of filters that satisfy the Boolean lifting property. The paper is organized into sections, each addressing distinct aspects of this relationship.

In Section 2, we introduce the concept of a pseudo-irreducible filter in a residuated lattice. After comparing this concept with other notions, such as prime filters, we provide several characterizations of such filters, which will be used in subsequent sections. In Sections 3 and 4, we revisit the definition of the Boolean lifting property for filters. Using pseudo-irreducible filters and the residuated lattice of fractions, we establish several useful characterizations for filters that satisfy the Boolean lifting property. Section 5 focuses on the Boolean lifting property of the radical of a filter. By introducing weak MTL-algebras and comparing them with other well-known classes of residuated lattices, we address an open question regarding the Boolean lifting property of the radical of a residuated lattice. Finally, we provide a topological solution to an open question about the Boolean lifting property of the radical of a residuated lattice [10, Open question 3.4].

\section{Preliminaries}  
In this section, we review some definitions and results which will be used throughout this paper.

\begin{definitie}\label{2.1}[7, 11, 21] A \textit{residuated lattice}  is an algebra $(L, \wedge, \vee, \odot, \rightarrow, 0, 1)$ 
	of type $(2, 2, 2, 2, 0, 0)$ satisfying the following axioms:\\
	(RL1) $(L, \wedge, \vee, 0, 1)$ is a bounded lattice (whose partial order is denoted by $\leq $);\\
	(RL2) $(L, \odot, 1)$ is a commutative monoid;\\
	(RL3) For every $x, y, z \in L$, $x \odot z \leq y$ if and only if $z\leq x\rightarrow y$ (residuation).
\end{definitie}
For   $x, y \in L$  and  $n\in \mathbb{N}$, we define \[x^{\ast}:= x\rightarrow 0,  x^{**}:=(x^*)^*, x^0:=1,  x^n:=x^{n-1}\odot x,   \text{  and }    x\leftrightarrow y:=(x\to y)\wedge(y\to x).\]

An element $x$  of a residuated lattice $L$ 
is called \textit{complemented} if there is an element $y \in L $ such that $x\wedge y = 0$ and $x\vee y= 1$, if such an element $y$ exists it is unique  and it is called  \textit{the  complement of} $x$. We  denote the  complement  of $x$ by $x'$.
The set of all complemented elements in $L$ is denoted by $ B(L)$ and
is called \textit{the Boolean center of}
$L$.

In the following proposition, we collect some main and well-known properties of residuated lattices, and throughout the paper, we frequently use them without referring.

\begin{propozitie}\label{2.2}  [7, 11, 21]
	Let $L$ be a residuated lattice, $x, y, z\in L$ and $e, f\in B(L)$. Then we have the following statements:\\
	\begin{enumerate}
		
		\item $x \leq y$ if and only if $ x \to y = 1$;
		\item   If $x \leq y$, then $y^*\leq x^*$;
		\item  $x\odot x^* = 0$;
		\item  $x\odot y=0$ if and only if $x\leq y^*$;
		\item  $x\leq x^{\ast\ast}$, $x^{\ast\ast\ast}=x^{\ast}$;
		\item   $x \to (y \to z) = y \to (x \to z) = (x\odot y) \to z$;
		\item $1\to x=x$ and $x\to 1=1$;
		\item $x\odot(x\to y)\leq y$;
		\item $x=x\leftrightarrow 1$ and $x^*= x\leftrightarrow 0$;
		\item  $(x\vee y)^* = x^* \wedge y^*$; 
		\item $x\vee(y\odot z)\geq (x\vee y)\odot (x\vee z)$;
		\item $x\odot(y\vee  z)= (x\odot y)\vee (x\odot z)$;
		\item $e' = e^*$,  $e\wedge e^*=e\odot e^*=0$, and $e^{**} = e$;
		\item $e\odot f= e \wedge f\in B(L)$;
		\item $e\odot x = e \wedge x$;
		\item $e \wedge (x \odot y) = (e \wedge x) \odot (e \wedge y)$;
		\item $e \vee (x \odot y) = (e \vee x) \odot (e \vee y)$;
		\item $e \wedge (x \vee y) = (e \wedge x) \vee (e \wedge y)$;
		\item $e \vee (x \wedge y) = (e \vee x) \wedge (e \vee y)$;
		\item $(x\wedge e)^*=x^*\vee e^*$;
		\item $e\to x=e^*\vee x$ and $x\to e=x^*\vee e$;
		\item $e\odot(e\to x)=e\wedge x$ and  $x\odot(x\to e)=e\wedge x$;
		\item For each $n\in \mathbb{N}$, $e^n=e$;
		\item If $x\vee x^*=1$, then $x\in B(L)$.
	\end{enumerate}
	
\end{propozitie}

\begin{definitie}\label{2.3} [2, 3, 14, 21] A residuated lattice $L$ is called\\
	\begin{enumerate}
		\item \textit{semi-G-algebra}, if	$(x^2)^{\ast} = x^{\ast}$
		for all	$x \in L$;
		\item 
		\textit{G}$\ddot{o}$\textit{del 
			algebra}  (\textit{G-algebra for short}), if	$x^{2} = x$	for all	$x\in L$;
		\item \textit{Involution}, if 	$x^{\ast\ast} = x$	for all	$x \in L$;
		\item \textit{MTL-algebra}, if $(x\rightarrow y)\vee (y\rightarrow x) = 1$ for all $x,y \in L$;
		
		\item \textit{De Morgan}, if $(x\wedge y)^*=x^*\vee y^*$ for all $x, y\in L$;
		\item \textit{Stonean}, if $x^{\ast}\vee x^{\ast\ast}=1$  for all $x\in L$.
	\end{enumerate}
\end{definitie}


\begin{definitie}\label{2.4}
	A non-empty subset $F$  of a residuated lattice $L$ is called a \textit{filter} 
	of $L$ if the following conditions hold:\\	
	(F1) If $x, y\in F$, then $x\odot y \in F$;\\
	(F2) If $x\leq y$ and $x\in F$, then $y \in F$.
\end{definitie}

\begin{definitie}
	A non-empty subset $F$  of a residuated lattice $L$ is called a \textit{deductive system} 
	of $L$ if the following conditions hold:\\	
	(D1) $1 \in F$;\\
	(D2) If $x, x\to y\in F$, then $y \in F$.
\end{definitie}
A non-empty subset $F$ of a residuated lattice $L$ is a filter if and only if it is a deductive system, see [21] for more details.

We denote by $   Filt(L)$ the set of all filters of $L$.
An filter $F$ is called \textit{proper} if $F \neq L$. Clearly a filter $F$ of $L$ is proper if and only if $0\not\in F$ if and only if for each $x\in L$ we have either $x\not\in F$ or $x^*\not\in F$.

For a non-empty subset $S$ of a residuated lattice  $L$, we set $[S):=\bigcap\{F\in    Filt(L)\mid S\subseteq F\}$, that is called \textit{the filter of $L$ generated by $S$}. We denote by $[x)$ \textit{the filter of $L$ generated by} $\{x\}$.
For $F\in    Filt(L)$ and $x\in L$, we set $F (x) := [F\cup \{x\})$. It is well-known that the  lattice $(   Filt(L),\subseteq)$ is distributive and complete. Actually,  for a family $\{F_i\}_{i\in A}$ of filters of $L$ we have $\bigwedge_{i\in A} F_i=\bigcap_{i\in A}F_i $ and $\bigvee_{i\in A} F_i=[\bigcup_{i\in A}F_i )$. If $F, G\in    Filt(L)$ we set $F\to G:=\{x\in L\mid F\cap [x)\subseteq G\}$ (see  [21, 24] and the references given there).

\begin{propozitie}\label{2.5}\label{3.7} [21] Let  $S$ be a non-empty subset of a residuated lattice $L$, $x, y \in L$ and $F, G\in    Filt(L)$. Then
	
	\begin{enumerate}
		\item $[S) = \{x\in L \mid  s_{1} \odot \cdots \odot s_{n}\leq x $ \text{  }\text{for some}\text{  } $n \geqslant 1 \text{  }\text{and}\text{  } s_{1},...,s_{n}\in S\}$;
		\item $[x) = \{z\in L \mid x^n\leq z $ for some $n \geqslant 1\}$. In particular, if $e\in  B(L)$, then $[e) = \{z\in L \mid e\leq z \}$;
		\item $F (x) = \{z\in L \mid  i\odot x^n\leq z $ for some $i \in I  $ and $  n\geqslant 1\}$;
		\item $F (x)\vee F (y)=F(x\wedge y)=F(x\odot y)$. In particular, $[x)\vee [y)= [x \wedge y) = [x \odot y)$;
		\item $F (x)\cap F (y)=F(x\vee y)$. In particular, $[x)\cap [y) = [x \vee y)$;
		\item $F \vee G= [F \cup G) =\{x \in L \mid  a\odot b\leq x$ for some  $a\in  F$ and $b\in G\}$;
		\item $F\to G\in   Filt(L)$.
	\end{enumerate}
\end{propozitie}


If $F$ is a filter of a residuated lattice $L$, then the
binary relation $\theta_{F}$ on $L$ defined by $(x, y)\in\theta_{F}$  if and only if  $x\leftrightarrow y\in F$
is a congruence on $L$. If 
$L/F$ denotes the quotient set $L/\theta_{F}$, then $L/F$ becomes a residuated lattice with the natural
operations induced from those of $L$.
Also, for $x\in L$, we denote by $x/{F}$ the  class of $x$ concerning to $\theta_{F}$.


For a non-empty subset $X$ of a  residuated lattice $L$, we set $X/F :=\{x/F\mid x\in X\} $.  Clearly, for $x \in L$;
$x/F =1/F $ if and only if $ x \in F$, and $x/F= 0/F$ if and only if $ x^* \in F$. Also, $x/F\leq y/F$ if and only if $x\to y\in F$ for each $x, y\in L$. Furthermore, if $G$ is a filter of $L$ containing $F$, then $G/F$ is a filter of $L/F$ and for $x\in L$ we have $x/F\in G/F$ if and only if $x\in G$.

Recall  that a proper filter $P \in    Filt(L)$ is called \textit{prime}
if  for $x, y\in L$, $x\vee y\in P$ implies either $x\in P$ or $y \in P$. We denote by $  SpecF(L)$ the set of all prime filters of $L$. An easy argument shows that a proper filter $P$ is prime if and only if  for $F, G\in    Filt(L)$ if $F\cap G\subseteq P$ then we have either $F\subseteq P$ or $G\subseteq P$. 	Also, a proper filter $M\in    Filt(L)$ is called \textit{maximal} if $M$ is not strictly contained in a proper filter
of $L$. We denote by $    MaxF(L)$ the set of all maximal filters of $L$. Clearly, every maximal filter  of  a 
residuated lattice is prime. Also, every proper filter is contained in a maximal filter, see  [11, 21, 24] for more details. 

For every subset $X$ of a residuated lattice $L$, we set $V (X) := \{P \in   SpecF(L) \mid X \subseteq P\}$, and for each $x\in L$, we set $V (x):=V (\{x\})$.  The family $\{V (X)\}_{X\subseteq L }$  satisfies the axioms for closed sets for
a topology  over $  SpecF(L)$. This topology is called  \textit{the Stone topology}. Since every maximal filter  is prime, we can consider $    MaxF(L)$ as a subspace of $  SpecF(L)$. Actually, for each $X\subseteq L $, we set $V_{Max}(X):=V(X)\cap     MaxF(L)$, then the  family $\{V_{Max} (X)\}_{X\subseteq L }$ satisfies the axioms for closed sets for
a topology over $  MaxF(L)$, for more information see [8].

\begin{propozitie}\label{2.70}[21] 
	Let $L$ be a  residuated lattice, $F\in    Filt(L)$ and $a\in L\setminus F$. Then we have the following statements:
	\begin{enumerate}
		\item There is a prime filter $P$ of $L$ such that $F\subseteq P$ and $a\not\in P$;
		\item $F$ is the intersection of all prime filters which contain $F$, that is,  $F=\bigcap V(F)$.
		
	\end{enumerate}
	
\end{propozitie}

\begin{propozitie}\label{sm}
	Let $L$ be a residuated lattice, $F, G\in   Filt(L)$ and $x, y\in L$. Then the following statements hold:
	\begin{enumerate}
		\item $F = L$ if and only if $V_{Max}(F)=\varnothing$ if and only if $V(F)=\varnothing$;
		\item $V(F) \cup V(G)=V(F \cap G)$ and $V_{Max}(F) \cup V_{Max}(G)=V_{Max}(F \cap G)$;
		\item  $V(x) \cup V(y) =V(x\vee y)$ and  $V_{Max}(x) \cup V_{Max}(y) =V_{Max}(x\vee y)$;
		\item If $\{X_k\}_{k\in K}$ is a family of subsets of $L$, then $V(\bigcup_{k\in K}X_k)=\bigcap_{k\in K}V(X_k)$ and   $V_{Max}(\bigcup_{k\in K}X_k)=\bigcap_{k\in K}V_{Max}(X_k)$.
	\end{enumerate}
\end{propozitie}

	\section {Pseudo-irreducible filters in residuated lattices}\label{1} 
	
	In this section, we first introduce the notion of  pseudo-irreducible filters in a residuated lattice and after that, we consider some of their properties.

	\begin{definitie}\label{3.1}
		A \textit{ pseudo-irreducible filter}  of a residuated lattice  $L$ is a proper filter $F$ of $L$ such that  for each $G, H\in    Filt(L)$, if $F=G\cap H$ and $G\vee H=L$, then we have either $G=L$ or $H=L$.
	\end{definitie}
	
	
	\begin{exemplu}\label{3.2}Let $L = \{0, a, b, c,  1\}$ be a lattice whose  Hasse diagram is given in the following  figure:
		
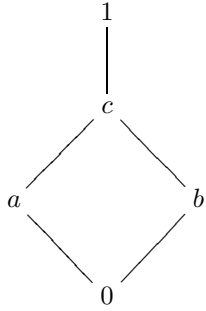
\[
		\xymatrix{
			& & 1\ar@{-}[d]\\
			&   &c \ar@{-}[rd]\ar@{-}[ld] &    \\ 
			&a \ar@{-}[rd]& & b \ar@{-}[ld]\\
			&& 0}		\]
		Consider the operations  $\odot$ and  $\rightarrow$ given by the
		following tables, see  [15, Example 4.4]:		
		\begin{center} 
			\begin{tabular}{r|l l l l l }
				$\odot$ & 0 & $a$ & $b$ &$c$  &  1 \\\hline
				0& $0$ &$0$& $0$&$0$&$0$\\
				$a$ &$0$ &$a$ & $0$& $a$&$a$\\
				$b$ &$0$ &$0$& $b$&$b$&$b$\\
				$c$ & $0$ &$a$& $b$&$c$&$c$\\
				1 & $0$& $a$&$b$& $c$& $1$\\
			\end{tabular} 
			\hspace{1cm}
			\begin{tabular}{r|l l l l l }
				$\rightarrow$ & 0 & $a$ & $b$ &$c$  & 1 \\\hline
				0& $1$ &$1$& $1$&$1$&$1
				$\\
				$a$ &$b$ &$1$ & $b$& $1$&$1$\\
				$b$ &$a$ &$a$& $1$&$1$&$1$\\
				$c$ & $0$ &$a$& $b$&$1$&$1$\\
				1 & $0$& $a$&$b$& $c$& $1$\\
			\end{tabular}
			
			\vspace{0.5cm}
		\end{center}
		\noindent	 Set $G:=[a)=\{a, c, 1\}$  and $H:=[b)=\{b, c, 1\}$. Since $[c)=\{c, 1\}=G\cap H$ and $G\vee H=L$, we have $[c) $ is not a pseudo-irreducible filter of $L$.
	\end{exemplu}
	Recall that a residuated lattice $L$ is called \textit{local} if $L$ has a unique maximal filter.
	
	In the following proposition, we consider conditions under which all proper filters of a residuated lattice are  pseudo-irreducible.
	\begin{propozitie}\label{om}
		Every proper filter of a  residuated lattice $L$ is  pseudo-irreducible if and only if $L$ is local.
	\end{propozitie}
	
	\begin{proof} 
		
		$\Rightarrow).$	 Assume that every proper filter of  $L$ 
		is pseudo-irreducible. If $M$ and $N$ are two distinct maximal filters of $L$, then $F:=M\cap N$ is pseudo-irreducible by our assumption. Now since $M\vee N=L$, we deduce that either $M=L$ or $N=L$, which is impossible. Hence,  $L$ has a unique maximal filter.
		
		$\Leftarrow).$ Let $L$ be a local residuated lattice with the unique maximal filter  $M$. If  $G, H\in    Filt(L)$ such that 
		$G\vee H=L$, then we have either $G\nsubseteq M$ or $H\nsubseteq M $, or equivalently, we have either  $G=L$ or $H=L$.  Therefore, every proper filter of  $L$ 
		is pseudo-irreducible. 		
	\end{proof}

	
	\begin{propozitie}\label{pmi} 
		Every prime ( or maximal)  filter of a residuated lattice is   pseudo-irreducible.
	\end{propozitie}
	\begin{proof} Let $F$ be a prime filter of  a residuated lattice ${L}$. Then $F$ is a proper filter. If there are $G, H\in    Filt(L)$ with
		$F=G\cap H$ and $G\vee H=L$, then  we have either $G\subseteq F$ or $H\subseteq F$. Assume that $G\subseteq F$. Hence, $G\subseteq H$, and we have $H=G\vee H=L$, that is, $H=L$.  Similarly, if  $H\subseteq F$, then we have $G=L$. Therefore,  $F$  is 
		pseudo-irreducible. Finally, since every maximal filter is prime, we deduced that every maximal filter is also pseudo-irreducible.
	\end{proof}
	
	In the following example, we show that  pseudo-irreducible filters
	need not be prime in general. 
	\begin{exemplu}\label{dd} Let $L = \{0, a, b, n, c, d, 1\}$ be a lattice whose  Hasse diagram is given in the following figure:
		
		\begin{center}
			$
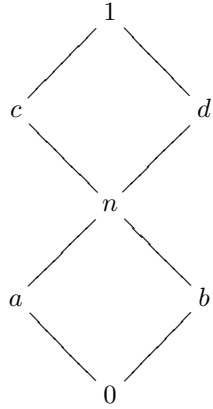
\xymatrix{
				& & 1\ar@{-}[rd]\ar@{-}[ld]\\
				& c \ar@{-}[rd] & & d\ar@{-}[ld]  & \\
				& &n\ar@{-}[rd] \ar@{-}[ld]&&\\
				&a \ar@{-}[rd]  & &b \ar@{-}[ld] &\\
				&& 0}$
		\end{center}
		
		Consider the operations  $\odot$ and  $\rightarrow$ given by the
		following tables (see [17, Example, Page 190]):
		
		\begin{center}\begin{tabular}{r|l l l l l l l }
				$\odot$ & 0 & $a$ & $b$ &$n$  & $c$ & $d$ & $1$   \\\hline
				0 & $0$ & $0$ &$0$  & $0$ & $0$ & $0$  & $0$ \\
				$a$& $0$ & $a$ &$0$  & $a$ & $a$ & $a$  & $a$\\
				$b$ & $0$ & $0$ &$b$  & $b$ & $b$ & $b$  & $b$\\
				$n$ & $0$ & $a$ &$b$  & $n$ & $n$ & $n$  & $n$ \\
				$c$ & $0$ & $a$ &$b$  & $n$ & $c$ & $n$  & $c$\\
				$d$ & $0$ & $a$ &$b$  & $c$ & $c$ & $d$  & $d$ \\
				1 & 0 &$a$ & $b$ &$n$  & $c$ & $d$ & $1$ 
			\end{tabular}
			\vspace*{0.5cm}
			\begin{tabular}{r|l l l l l l l }
				$\rightarrow$  &0& $a$ & $b$ &$n$  & $c$ & $d$ & $1$   \\\hline
				0 & $1$ & $1$ &$1$  & $1$ & $1$ & $1$  & $1$ \\
				$a$& $b$ & $1$ &$b$  & $1$ & $1$ & $1$  & $1$\\
				$b$ & $a$ & $a$ &$1$  & $1$ & $1$ & $1$  & $1$\\
				$n$ & $0$ & $a$ &$b$  & $1$ & $1$ & $1$  & $1$ \\
				$c$ & $0$ & $a$ &$b$  & $d$ & $1$ & $d$  & $1$\\
				$d$ & $0$ & $a$ &$b$  & $n$ & $n$ & $d$  & $1$ \\
				1 & 0 &$a$ & $b$ &$n$  & $c$ & $d$ & $1$ 
			\end{tabular}

			\vspace{0.5cm}
		\end{center}
		Since $c\vee d=1\in\{1\}$ and $c, d\not\in\{1\}$, we conclude  that $\{1\}$ is not a prime filter of $L$.	Now if $G, H\in    Filt(L)$ such that $G\vee H=L$, then by Proposition 1.6, there are $x\in G$ and $ y\in H$ such that $x\odot y=0$. If $G$ and $H$ are proper filters, then by the above tables, we have $x\not= y$ and $x, y\in \{a, b\}$. It follows that either ($a\in G$ and $b\in H$)  or ($b\in G$ and $a\in H$). Now since $a\vee b=n$, in both cases, we have $n\in G\cap H$.
		Therefore, there are no proper filters $G, H\in    Filt(L)$ such that $\{1\}=G\cap H$ and $G\vee H=L$, and so the filter $\{1\}$ is pseudo-irreducible.
	\end{exemplu}	
	
	We recall that a residuated lattice $L$ is \textit{directly indecomposable} if $L\cong L_{1}\times L_{2 }$ implies either
	$L_{1}$ or $L_{2 }$ is trivial, where $L_{1}$ and  $L_{2 }$ are two residuated lattices and  $L_{1}\times L_{2 }$ is their
	direct product.  It is well-known that a non-trivial residuated lattice $L$  is  directly indecomposable if and only if $ B(L) = \{0, 1\}$. For more details, we refer the reader to [4] and [7]. 
	\begin{propozitie}[8, Proposition 6.8]\label{ssss}
		For a residuated lattice $L$, we have \[   Clop(  SpecF(L))=\{V(e)\mid e\in  B(L)\},\] where $   Clop(  SpecF(L))$ is the set of all closed and open subsets of $  SpecF(L)$ with respect to the Stone topology. 
	\end{propozitie}
	
	\begin{propozitie}\label{ttt}
		Let $L$ be a residuated lattice. Then $L$ is directly indecomposable if and only if $  SpecF(L)$ is a connected topological space  with respect to the Stone topology.
	\end{propozitie}
	\begin{proof}
		Let $e, f\in  B(L)$. By  Proposition 1.6, we have  $e=f$ if and only if $[e)=[f)$. Also by  Proposition 1.7, $[e)=[f)$ if and only if $V(e)=V(f)$. Consequently,  $e=f$ if and only $V(e)=V(f)$.
		
		On the other hand, we know that $L$ is directly indecomposable if and only if $ B(L) = \{0, 1\}$. Now from Proposition 2.6 and the fact that a topological space is connected if and only if it has no non-trivial clopen subsets, we conclude that $L$ is directly indecomposable if and only if $ B(L)=\{0, 1\}$ if and only if $   Clop(  SpecF(L))=\{V(0)=\varnothing, V(1)=  SpecF(L)\}$ if and only if $  SpecF(L)$ is a connected topological space  with respect to the Stone topology.  
	\end{proof}	
	
	In the following theorem, we provide some conditions that are equivalent to a filter being pseudo-irreducible.
	\begin{teorema}\label{eqm} 
		The following statements are equivalent for a proper filter $F$ of a  residuated lattice $L$:
		\begin{enumerate}
			\item  $F$ is  a pseudo-irreducible filter of $L$; 
			
			\item For $x, y\in L $, if $x\vee y\in F$ and $x\odot y=0$, then we have  either $x\in F$ or $y\in F$;
			\item For $x\in L $, if $x\vee x^*\in F$, then we have  either $x\in F$ or $x^*\in F$;
			
			\item The  residuated lattice $L/F$ is  non-trivial and directly indecomposable;
			\item $V(I)$ is connected as a subspace of $  SpecF(L)$ with respect to the Stone topology.
			
		\end{enumerate}
		
	\end{teorema}
	\begin{proof}
		$(1)\Rightarrow (2)$. Assume that $x\vee y\in F$ and $x\odot y=0$ for some $x, y\in L$.  By Proposition 1.6,
		we have $F=F(x\vee y)=F (x)\cap F (y)$ and
		$F(x)\vee F(y)=   F(x\odot y)=F(0)=L$. By assumption, we conclude that either $F(x)=L$ or $F(x)=L$. If $F(x)=L$, then $y\in L=F(x)$. Since $y\in F(y)$, we have $y\in F (x)\cap F (y)=F$.
		Similarly, if $F(y)=L$, then we have $x\in F$.
		
		$(2)\Rightarrow(1)$. Suppose that  $G, H\in    Filt(L)$ such that $F=G\cap H$ and $G\vee H=L$.
		Since $0\in L=G\cap H$, we have $a\odot b=0$ for some $a\in G$ and $b\in  H$ by Proposition 1.6.
		So $a\vee b\in G\cap H=F$. Our assumption implies that  either $a\in F$ or $b\in F$. Suppose that
		$a\in F$. Since $F\subseteq H$, we have $a\in H$. Thus we conclude that $a\odot b\in H$. Therefore, $0\in H$ and so $H=L$. Similarly, if $b\in F$, then $G=L$. Therefore, $F$ is  a pseudo-irreducible filter of $L$.

		$(2)\Rightarrow(3)$. It is clear since for each $x\in L$, we have $x\odot x^*=0$.
		
		$(3)\Rightarrow(2)$. Assume that $x\vee y\in F$ and $x\odot y=0$ for some $x, y\in L$. Then, $x\leq y^*$ and $y\leq x^*$. It follows that $(x\vee y)\leq (x\vee x^*)$ and $(x\vee y)\leq (y\vee y^*)$, and consequently  $x\vee x^*, y\vee y^*\in F$. By assumption, we have  (either $x\in F$ or $x^*\in F$) and  (either $y\in F$ or $y^*\in F$). If $x^*, y^*\in F$, then $(x\vee y)^*=x^*\wedge y^*\in F$. Hence, $x\vee y, (x\vee y)\to 0\in F$. This implies $0\in F$, a contradiction. Therefore, we have either  $x\in F$ or $y\in F$.

		$(3)\Rightarrow (4)$.	Since $F$ is a proper filter of $L$, the  residuated lattice $L/F$ is  non-trivial. Now we want to prove that ${L}/{F}$ is directly indecomposable. For this purpose, we show that 
		$ B({L}/{F}) =\{{0}/{F}, {1}/{F}\}$.
		Let ${x}/{F}\in  B({L/}{F})$. Then $ {x}/{F}\vee{x^{\ast}}/{F}={1}/{F}$. Hence, ${(x\vee x^{\ast})}/{F}={1}/{F}$, and so $x\vee x^{\ast}\in F$. By our assumption, we have  either $x\in F$ or $x^*\in F$, or equivalently, we have  either $x/F=1/F$ or $x/F=0/F$. It follows that 
		$ B({L}/{F}) =\{{0}/{F}, {1}/{F}\}$.

		$(4)\Rightarrow (3)$. Suppose $x\vee x^*\in F$ for some $x\in L$. Thus $x/F\vee x^*/F=(x\vee x^*)/F=1/F$, and so $x/F\in  B(L/F)$ by Proposition 1.2. Now since $L/F$ is directly indecomposable, we conclude that either $x/F=1/F$ or $x/F=0/F$, or equivalently, we have either $x\in F$ or $x^*\in F$. 
		
		$(4)\Leftrightarrow (5)$. We know the prime filters of ${L}/{F}$
		are exactly of the form ${P}/{F}(:=\{{x}/{F}\mid x\in P\})$, where $P$ is a prime filter of $L$ containing $F$. Hence, the space $V(F)$ is homeomorphic to the space $  SpecF({L}/{F})$ by the natural homeomorphism. Using  Proposition 2.7 and the above argument, we have   ${L}/{F}$ is directly indecomposable if and only if  $  SpecF({L}/{F})$ is a connected topological space with respect to the Stone topology if and only if  $V(F)$ is connected as a subspace of $  SpecF(L)$ with respect to the Stone topology.
	\end{proof}

	The following corollary is a direct consequence of Theorem 2.8.
	\begin{corolar}\label{pib}
		If $F$ is a pseudo-irreducible (e.g., maximal or prime by Proposition 2.4) filter of $L$, then for each $e\in B(L)$ we have either $e\in F$ or $e^*\in F$.
	\end{corolar}

	We end this section with the following proposition, which considers when the intersection of two pseudo-irreducible filters is pseudo-irreducible.
	\begin{propozitie}\label{ipi}
		Let $F$ and $G$ be two pseudo-irreducible filters of a  residuated lattice  $L$. Then 
		$G\vee H\neq L$ if and only if $G\cap H$ is  pseudo-irreducible.
		
		\begin{proof} $\Rightarrow)$. Let $x\vee x^*\in G\cap H$. Then $x\vee x^*\in G$ and $x\vee x^*\in  H$. Thus by Theorem 2.8, we have (either $x\in G$ or $x^*\in G$) and (either $x\in H$ or $x^*\in H$). If $x\in G$ and $x^*\in H$, then $0=x\odot x^*\in G\vee H$. It follows that $G\vee H= L$,  which is impossible. Similarly, if $x^*\in G$ and $x\in H$, then  we have $G\vee H= L$, which is impossible.  Hence, we have either $x\in G\cap H$ or $x^*\in G\cap H$. Therefore, $G\cap H$ is  pseudo-irreducible by Theorem 2.8
			
			$\Leftarrow)$. Assume that $G\cap H$ is  pseudo-irreducible. If $G\vee H= L$, then we have either $G=L$ or $H=L$. It is impossible since $F$ and $G$ are pseudo-irreducible filters of $L$ and they must be proper.
		\end{proof}
	\end{propozitie}
	\section{Boolean Lifting Property and pseudo-irreducibility} \label{2}  
	In this section, we consider the Boolean lifting property for a filter and its  connection with pseudo-irreducibility. We begin with the following definition.
	
	\begin{definitie}\label{lll} A (not necessarily proper) filter $F$ of a residuated lattice $L$ 
		has \textit{the Boolean lifting property (BLP)}, whenever each $\alpha\in  B(L/F )$ can be lifted, that is,  there exists 
		$e \in  B(L)$ such that $\alpha=e/F$.
	\end{definitie}
	\begin{exemplu}\label{eg}
		An easy argument shows that the  improper filter and the trivial filter of a residuated lattice have BLP, see also [8, Corollary 4.5].
	\end{exemplu}
	
	First of all, in the following theorem, we provide some conditions that are equivalent to a filter having BLP. Later, we will present additional equivalent conditions.
	\begin{teorema}\label{a1}
		Let $F$ be a filter of a residuated lattice $L$. Then the following statements are equivalent:
		\begin{enumerate}
			\item $F$ has  BLP;
			
			\item For each $x\in L$, if $x\vee x^*\in F$, then there exists $e\in  B(L)$ such that $x\leftrightarrow e\in F$,  or equivalently, $x\vee e^*, x^*\vee e\in F$. In this case, $x^*\leftrightarrow e^*\in F$; 
			\item If there are $G, H\in   Filt(L)$ such that $F=G\cap H$ and $G\vee H=L$, then there exists $e\in  B(L)$ such that $e\in G$ and $e^*\in H$;
			\item  If there are $G, H\in   Filt(L)$ such that $G\cap H\subseteq F$ and $G\vee H=L$, then there exists $e\in  B(L)$ such that $e\in F\vee G$ and $e^*\in F\vee H$.
		\end{enumerate}
	\end{teorema} 
	\begin{proof}  
		$(1)\Rightarrow (2)$. Assume that $x\vee x^*\in F$. Thus $x/F\vee x^*/F=(x\vee x^*)/F=1/F$, and so $x/F\in  B(L/F)$ by Proposition 1.2. By our assumption, there exists 
		$e \in  B(L)$ such that $x/F=e/F$. This implies that 	$x\leftrightarrow e\in F$, or equivalently, $x\vee e^*, x^*\vee e\in F$. If this case happens, since $x/F=e/F$, we have $x^*/F=e^*/F$, or equivalently, 	$x^*\leftrightarrow e^*\in F$.
		
		$(2)\Rightarrow (1)$.	If $x/F\in  B(L/F)$, then $(x\vee x^*)/F=x/F\vee x^*/F=1/F$, and so $x\vee x^*\in F$. Hence by our assumption, there exists $e\in  B(L)$ such that $x\leftrightarrow e\in F$. It follows that $x/F=e/F$, and so $F$ has  BLP.
		
		$(2)\Rightarrow (3)$. Assume that there are $G, H\in   Filt(L)$ such that $F=G\cap H$ and $G\vee H=L$. Since $G\vee H=L$, $0\in G\vee H$. Thus by Proposition 1.6, there exist $x\in G$ and $y\in H$ such that $x\odot y=0$ and $x\vee y\in G\cap H=F$. We conclude from  $x\odot y=0$ that $x\leq y^*$ and $y\leq x^*$. This gives $x\vee y\leq x\vee x^*$, and consequently $x\vee x^*\in F$. By assumption, there exists $e\in  B(L)$ such that $x\leftrightarrow e\in F=G\cap H$ and $x^*\leftrightarrow e^*\in F=G\cap H$. From $x\leftrightarrow e\in G$ and $x\in G$, we have $e\in G$. Also,  since $y\leq x^*$, we have $x^*\in H $. Now since $x^*\leftrightarrow e^*\in H$, we conclude that $e^*\in H$. Hence, $e\in G$ and $e^*\in H$.
		
		$(3)\Rightarrow (2)$. Assume that $x\vee x^*\in F$ for some $x\in L$. Using Proposition 1.6, we have $F=F(x\vee x^*)=F(x)\cap F(x^*)$ and  $F(x)\vee F(x^*)=F(x\odot x^*)=F(0)=L$. By our assumption, there exists $e\in  B(L)$ such that $e\in F(x)$ and $e^*\in F(x^*)$. We conclude from $e\leq x\to e$ and $x\leq e\to x$ that $x\to e, e\to x\in F(x)$, or equivalently, $x\leftrightarrow e\in F(x)$. Also, from $x^*\leq x\to e$ and $e^*\leq e\to x$, we have $x\to e, e\to x\in F(x^*)$, or equivalently, $x\leftrightarrow e\in F(x^*)$. Therefore, $x\leftrightarrow e\in F(x)\cap F(x^*)=F$.
		
		$(3)\Leftrightarrow (4)$. It is clear from the fact that the lattice $(   Filt(L), \subseteq)$ is distributive.
	\end{proof}

	A quick consequence of Theorems 2.8 and 3.3, is the following corollary.
	\begin{corolar}\label{psblp}
		Every pseudo-irreducible  filter has  BLP. 
	\end{corolar}

	The converse of the above corollary is not true in general.
	\begin{exemplu}\label{a2}
		Consider any directly decomposable residuated lattice $L$ (that is, $L=L_1\times L_2$, where $L_1$ and $L_2$ are two non-trivial residuated latices). Then by Theorem 2.8, the filter $\{(1, 1)\}$ is not pseudo-irreducible, but it has BLP by Example 3.2.   
	\end{exemplu}
	In the following proposition, we consider conditions for a residuated lattice under which each proper filter that has BLP, is a pseudo-irreducible filter.
	\begin{propozitie}\label{a3}
		Let $L$ be a non-trivial residuated lattice. Then $L$ is directly indecomposable if and only if every proper filter  that  has BLP  is a pseudo-irreducible filter.
	\end{propozitie}
	\begin{proof} $\Rightarrow)$. Let $F$ be a proper filter of $L$  that  has BLP and  $x\vee x^*\in F$ for some $x\in L$. By assumption and Theorem 3.3, there exists $e\in  B(L)$ such that $x\leftrightarrow e\in F$, or equivalently, $x\vee e^*, x^*\vee e\in F$. Since $L$ is directly indecomposable, we have $ B(L)=\{0, 1\}$. It follows that either $e=0$ or $e=1$. Hence, we have either $x\in F$ or $x^*\in F$. Consequently, $F$ is a pseudo-irreducible filter of $L$ by Theorem 2.8.
		
		$\Leftarrow)$. By Example 3.2, the filter $\{1\}$ is a proper filter of $L$ that has BLP. Hence, it is pseudo-irreducible. Now since $L\cong L/\{1\}$,  we conclude that $L$ is directly indecomposable by Theorem 2.8.
	\end{proof}

	In the following propositions, we will consider some relationship between our two main concepts: pseudo-irreducibility of filters and having BLP, in this paper.
	
	\begin{propozitie}\label{b8}
		Let $F$ and $G$ be two filters of $L$ such that $F\vee G\not=L$. Then if
		$F$ is a pseudo-irreducible filter and $G$  has BLP, then $F\cap G$ has BLP.
	\end{propozitie}
	
	\begin{proof}
		Assume that $x\vee x^*\in F\cap G$ for some $x\in L$. Since $x\vee x^*\in G$,  there exists $e\in  B(L)$ such that $x\leftrightarrow e, x^*\leftrightarrow e^* \in G$ by Theorem 3.3. On the other hand, since $x\vee x^*\in F$, by Theorem 2.8  we have either $x\in F$ or $x^*\in F$.  Also, by using Corollary 2.9 and  pseudo-irreducibility of $F$, we have either $e\in F$ or $e^*\in F$. Now we consider the following cases:
		
		$\mathbf{Case\text{ } 1}$: If $x, e\in F$, then $x\leftrightarrow e\in F$ and so $x\leftrightarrow e\in F\cap G$.
		
		$\mathbf{Case\text{ } 2}$: If $x^*, e^*\in F$, then $x\leftrightarrow e\in F$ since $x^*\leq x\to e$ and $e^*\leq e\to x$. Thus $x\leftrightarrow e\in F\cap G$.
		
		$\mathbf{Case\text{ } 3}$: If $x, e^*\in F$, then $x\odot (x\to e)\in F\vee G$. From $x\odot (x\to e)\leq e$, we conclude that $e\in F\vee G$. Also, since $e^*\in F\subseteq F\vee G$, we have $0=e\odot e^*\in F\vee G$, or equivalently, $F\vee G=L$, which is a impossible.
		
		$\mathbf{Case\text{ } 4}$:  If $x^*, e\in F$, then $x^*\odot (x^*\to e^*)\in F\vee G$. From $x^*\odot (x^*\to e^*)\leq e^*$, we deduce that $e^*\in F\vee G$. Also, since $e\in F\subseteq F\vee G$, we have $0=e\odot e^*\in F\vee G$, or equivalently, $F\vee G=L$, which is a impossible.
		
		Therefore, $F\cap G$ has BLP by Theorem 3.3.
	\end{proof}
	
	\begin{propozitie}\label{b9}
		If $F$ and $G$ are two proper filters of $L$ with $F\subseteq G$. Then we have the following statements:
		\begin{enumerate}
			\item  If $F$ is   pseudo-irreducible  and $G$ has BLP, then $G$ is  pseudo-irreducible; 
			
			\item  $G$ is a pseudo-irreducible filter of $L$ if and only if $G/F$ is a  pseudo-irreducible filter of $L/F$;
			
			\item If $G$  has BLP in $L$, then $G/F$  has BLP  in $L/F$. The converse is true if $F$ has BLP.	
		\end{enumerate}
	\end{propozitie}
	
	\begin{proof}
		$(1).$	Suppose that $x\vee x^*\in G$ for some $x\in L$. By assumption and Theorem 3.3, there exists $e\in  B(L)$ such that  $x\leftrightarrow e\in G$ and  $x^*\leftrightarrow e^*\in G$. Using Corollary 2.9, we have either $e\in F$ or $e^*\in F$. If $e\in F$, then $e\in G$. Hence $e, e\to x\in G$, and we have $x\in G$. Similarly, if $e^*\in F$, then $x^*\in G$. Therefore, $G$ is  a pseudo-irreducible filter of $L$ by Theorem 2.8. 
		
		$(2)$. Clearly $G$ is a proper filter of $L$ if and only if $G/F$ is a proper filter of $L/F$. By the second isomorphism theorem, we have $L/G\cong (L/F)/(G/F)$. Hence, by Theorem 2.8, $G$ is a pseudo-irreducible filter of $L$ if and only if the residuated lattice $L/F$ is directly indecomposable if and only if the residuated lattice $(L/F)/(G/F)$ is directly indecomposable if and only if $G/F$ is a  pseudo-irreducible filter of $L/F$.

		$(3)$.	Assume that $G$ has BLP in $L$ and $x/F\vee x^*/F\in G/F$ for some $x\in L$. Hence $(x\vee x^*)/F\in G/F$, and so $x\vee x^*\in G$. Now using Theorem 3.3,  there exists $e\in  B(L)$ such that $x\leftrightarrow e \in G$. Thus $(x\leftrightarrow e)/F\in G/F$. This shows that $x/F\leftrightarrow e/F\in G/F$ and $e/F\in  B(L/F)$. Therefore, $G/F$  has BLP  in $L/F$ by Theorem 3.3.
		
		Now suppose that $F$ has BLP in $L$ and $G/F$  has BLP  in $L/F$. If  $x\vee x^*\in G$, then $x/F\vee x^*/F=(x\vee x^*)/F\in G/F$. By assumption and Theorem 3.3, there is $\alpha \in  B(L/F)$ such that $x/F\leftrightarrow \alpha\in G/F$. Now since $\alpha \in  B(L/F)$ and  $F$ has BLP, there exists $e\in  B(L)$ such that $\alpha=e/F$. Thus we conclude that $x/F\leftrightarrow e/F\in G/F$, or equivalently, $(x\leftrightarrow e)/F\in G/F$. Hence $x\leftrightarrow e\in G$. This shows that $G$  has BLP in $L$ by Theorem 3.3.
	\end{proof}

	
	\section{Boolean Lifting Property and the residuated lattice of fractions}\label{3}
	First of all, in this section, we recall the definition of the residuated lattice of fractions of a residuated lattice  relative to a $\wedge-$closed system  and extend some of its theory from the viewpoint of filter theory. Then, we use the obtained results to provide some additional characterizations for filters that have BLP.
	
	Recall that a \textit{$\wedge-$closed system} of a residuated lattice $L$ is a non-empty subset $S$ of $L$ such that $1\in S$ and if
	$ x, y \in S$, then $x \wedge y \in S$. For a $\wedge-$closed system $S$ of a residuated lattice $L$,  we consider the binary relation $\theta_S$ defined by $(x, y)\in \theta_S$ if and only if  there
	exists $e\in S \cap  B(L)$ such that $x \wedge e = y \wedge e$. By [1, Lemma 4], $\theta_S$  is a congruence on $L$, and hence, $L[S]:=\frac{L}{\theta_S}$, the set of all equivalence class of $L$ with respect to $\theta_S$, becomes a residuated lattice with the natural
	operations induced from those of $L$. $L[S]$ is called \textit{the residuated lattice of fractions of $L$  relative to  $S$}. For each $x\in L$, $\frac{x}{S}$  denotes the equivalence class of $x$ relative to $\theta_S$, see [1] for more details.
	
	\begin{remarca}
		In the following, to avoid any misunderstanding between the symbols used for the two concepts \lq\lq the residuated lattice of fractions relative to a $\wedge-$closed system" and \lq\lq the quotient of a residuated lattice with respect to a filter", for a residuated lattice $L$,  we will use the symbol $L[S]$ (and $\frac{x}{S}$ for its equivalence class) for the residuated lattice of fractions of $L$ relative to a $\wedge-$closed system  $S$, and we will use the symbol $L/F$ (and $x/F$ for its equivalence class) for the quotient of the residuated lattice $L$ with respect to a filter $F$.
	\end{remarca}

	\begin{propozitie}\label{wedg}[1, Remark 6 and Proposition 4]
		For a $\wedge-$closed system $S$ of a residuated lattice $L$ we have the following statements:
		\begin{enumerate}
			\item $\frac{e}{S}=\frac{1}{S}$ for each $e\in S\cap B(L)$;
			\item $\frac{x}{S}\in B(L[S])$ 	if and only if $e\leq x\vee x^*$ for some
			$e\in S\cap  B(L)$.
			
		\end{enumerate}
	\end{propozitie}
	In the following proposition, we consider the filters of the residuated lattice of fractions relative to a $\wedge-$closed system.
	
	\begin{propozitie}\label{b1}
		Let $F$ be a filter  and $S$ be a $\wedge-$closed system  of a residuated lattice $L$. Then $F[S]:=\{\frac{x}{S}\mid x\in F\}$ is a filter of $L[S]$. Actually, every filter of $L[S]$ is of the form $F[S]$ for some filter $F$ of $L$.
	\end{propozitie}
	\begin{proof}
		Clearly, $F[S]\not=\varnothing$. Let $\frac{x}{S}, \frac{y}{S}\in F[S]$, where $x, y\in F$. Since $x\odot y\in F$, we have $\frac{x}{S}\odot\frac{y}{S}=\frac{x\odot y}{S}\in F[S]$. Now let $x\in F$  and $y\in L$ such that $\frac{x}{S}\leq\frac{y}{S}$. It follows that $\frac{x\to y}{S}=\frac{x}{S}\to\frac{y}{S}=\frac{1}{S}$. By definition, there exists $e\in S\cap  B(L)$ such that $(x\to y)\wedge e=1\wedge e=e$. So $e\leq x\to y$, or equivalently, $e\to(x\to y)=1$. From Proposition 1.2, we have $x\to(e\to y)=1\in F$, and so $e\to y\in F$ since $x\in F$. As $e\wedge(e\to y)=e\odot(e\to y)=e\wedge y$, we have $\frac{y}{S}=\frac{e\to y}{S}\in F[S]$. Hence,   $\frac{y}{S}\in F[S]$. Therefore, $F[S]$ is a filter of $L[S]$.
		
		Now assume that $T$ is a filter of $L[S]$. Set $F:=\{x\in L\mid \frac{x}{S}\in T\}$. It is easily seen that $F$ is a filter of $L$ and $T=F[S]$.
	\end{proof}

	A simple argument gives the following lemma. 
	\begin{lema}\label{b7}
		Let $F$ be a filter  and $S$ be a $\wedge-$closed system of a residuated lattice $L$. Then $\frac{x}{S}\in F[S]$ if and only if there exist $e\in S\cap  B(L)$ and $a\in F$ such that $x\wedge e=a\wedge e$, or equivalently, $x\vee e^*\in F$. In particular, $\frac{x}{S}=\frac{0}{S}$ if and only if there exists $e\in S\cap  B(L)$  such that $x\leq e^*$.
	\end{lema}
	
	\begin{propozitie}\label{b2}
		Let $F$ be a filter  and $S$ be a $\wedge-$closed system of a residuated lattice $L$. Then $F[S]=L[S]$ if and only if there exists $e\in S\cap  B(L)$ such that $e^*\in F$.
	\end{propozitie}
	
	\begin{proof}
		If $F[S]=L[S]$, then $\frac{0}{S}\in F[S]$. Hence $\frac{0}{S}=\frac{x}{S}$ for some $x\in F$, or equivalently, there is $e\in S\cap  B(L)$ such that $x\wedge e=0\wedge e=0$. It follows that $x\odot e=x\wedge e=0$, and so $x\leq e^*$. Now since $x\in F$, we have $e^*\in F$. 
		
		Conversely, if $e^*\in F$ for some $e\in S\cap  B(L)$, then $e^*\wedge e=0\wedge e$. Thus $\frac{0}{S}=\frac{e^*}{S}\in F[S]$. Therefore, $F[S]=L[S]$.
	\end{proof}
	
	\begin{propozitie}\label{b3}
		Let $F$ and $G$ be two filters  and $S$ be a $\wedge-$closed system of a residuated lattice $L$. Then $F[S]\cap G[S]=(F\cap G)[S]$. 
	\end{propozitie}
	\begin{proof}
		Clearly, $(F\cap G)[S]\subseteq F[S]\cap G[S]$. If $\frac{x}{S}\in F[S]\cap G[S]$, then by Lemma 4.4 there are $a\in F$, $b\in G$ and $e,f \in S\cap  B(L)$ such that $x\wedge e=a\wedge e$ and $x\wedge f=b\wedge f$. Hence we have
		\begin{align*}
			(a\vee b)	\wedge(e\wedge f)&=(a\vee b)\odot(e\odot f)\\&=(a\odot e\odot f)\vee (b\odot e\odot f)\\&=(a\wedge e\wedge f)\vee (b\wedge e\wedge f)\\&=(x\wedge e\wedge f)\vee(x\wedge e\wedge f)\\&=(x\wedge e\wedge f)\\&=x\wedge (e\wedge f).
		\end{align*}
		It follows that $\frac{x}{S}\in (F\cap G)[S]$ from the fact that $a\vee b\in F\cap G$ and $e\wedge f\in S\cap B(L)$. 
	\end{proof}

	\begin{propozitie}\label{b5}
		Let $F$  be a proper filter  and $S$ be a $\wedge-$closed system of a residuated lattice $L$ such that $F[S]\not=L[S]$. Then the following statements are equivalent:
		
		\begin{enumerate}
			\item $F[S]$ is a pseudo-irreducible filter of $L[S]$;
			\item If $e^*\vee x\vee x^*\in F$ for some $e\in S\cap B(L)$ and $x\in L$, then there exists $f\in S\cap B(L)$ such that we have either $f^*\vee x\in F$ or $f^*\vee x^*\in F$. 
		\end{enumerate}
	\end{propozitie}
	\begin{proof}
		By definition and Lemma 4.4, we first establish some facts about elements of residuated lattices $L$ and $ L[S]/F[S]$.
		\begin{align*}
			\frac{x}{S}/F[S]\in B(L[S]/F[S])&\Leftrightarrow \frac{x\vee x^*}{S}/F[S]=\frac{1}{S}/F[S]\\&\Leftrightarrow \frac{(x\vee x^*)}{S}\in F[S]\\& \Leftrightarrow (x\vee x^*)\wedge e=a\wedge e \text{  for some }  e\in S\cap  B(L) \text{ and }  a\in F\\& \Leftrightarrow e^*\vee x\vee x^*\in F \text{  for some }  e\in S\cap  B(L).
		\end{align*}
		Also, we have
		\begin{align*}
			\frac{x}{S}/F[S]=\frac{0}{S}/F[S]& \Leftrightarrow \frac{x^*}{S}\in F[S]  \\& \Leftrightarrow x^*\wedge f=a\wedge f \text{  for some }  f\in S\cap  B(L) \text{ and }  a\in F\\& \Leftrightarrow x^*\vee f^*\in F \text{  for some }  f\in S\cap  B(L),
		\end{align*}
		and
		\begin{align*}
			\frac{x}{S}/F[S]=\frac{1}{S}/F[S]&\Leftrightarrow \frac{x}{S}\in F[S] \\& \Leftrightarrow  x\wedge f=a\wedge f \text{  for some }  f\in S\cap  B(L) \text{ and }  a\in F\\& \Leftrightarrow x\vee f^*\in F \text{  for some }  f\in S\cap  B(L).
		\end{align*}
		Now by Theorems 2.8, we know that $F[S]$ is a pseudo-irreducible filter of $L[S]$ if and only if $L[S]/F[S]$ is directly indecomposable if and only if  $ B(L[S]/F[S])=\{\frac{0}{S}/F[S], \frac{1}{S}/F[S]\}$. Therefore,  the above facts complete the proof.
	\end{proof}

	\begin{propozitie}\label{b12}
		Let $F$ be a proper filter and $S$ be a $\wedge-$closed system of  a residuated lattice $L$. Then we have the following statements:
		\begin{enumerate}
			\item If $F$ is a  pseudo-irreducible filter of $L$ and $F[S]\not= L[S]$, then $F[S]$ is a   pseudo-irreducible filter  of $L[S]$;
			
			\item If $F$ has BLP in $L$, then $F[S]$  has BLP in $L[S]$.
		\end{enumerate}
		
	\end{propozitie}
	\begin{proof}
		$(1).$ Assume that $\frac{x}{S}\vee\frac{x^*}{S}\in F[S]$, hence we have  $\frac{x\vee x^*}{S}\in F[S] $. By Lemma 4.4, there exist $a\in F$ and $e\in S\cap B(L)$ such that $(x\vee x^*)\wedge e=a\wedge e$. Now Corollary 2.9 implies that either $e\in F$
		or  $e^*\in F$. Since  $F[S]\not=L[S]$, we have $e\in F$ by Proposition 4.5, hence  	$a\wedge e\in F$, and so $x\vee x^*\in F$ since $(a\wedge e)\leq (x\vee x^*)$. By assumption and Theorem 2.8, we have either $x\in F$ or $x^*\in F$, and hence we have either $\frac{x}{S}\in F[S]$ or $\frac{x^*}{S}\in F[S]$. Therefore, $F[S]$ is a   pseudo-irreducible filter  of $L[S]$.
		
		$(2).$ Assume that $\frac{x}{S}\vee\frac{x^*}{S}\in F[S]$. Hence $\frac{x\vee x^*}{S}\in F[S] $. By Lemma 4.4, there exist $a\in F$ and $e\in S\cap B(L)$ such that $(x\vee x^*)\wedge e=a\wedge e$. It follows that 
		\begin{align*}
			a\odot e=a\wedge e\leq x\vee x^*&\Rightarrow a\leq e\to(x\vee x^*)\\ &\Rightarrow a\leq e^*\vee (x\vee x^*).
		\end{align*}
		Hence, we have 
		\begin{align*}
			(x\vee e^*)\vee (x\vee e^*)^*&=(x\vee e^*)\vee (x^*\wedge e)\\&=(x\vee e^*)\vee (x^*\odot e)\\&\geq (x\vee e^*\vee x^*)\odot(x\vee e^*\vee e)\\&=(x\vee e^*\vee x^*)\odot 1\\&=(x\vee e^*\vee x^*)\geq a.
		\end{align*}
		Consequently, $a\leq (x\vee e^*)\vee (x\vee e^*)^* $, and since $a\in F$, we have $(x\vee e^*)\vee (x\vee e^*)^*\in F$. By hypothesis and Theorem 3.3, there exists $f\in  B(L)$ such that  $(x\vee e^*)\leftrightarrow f\in F$. Thus, $\frac{(x\vee e^*)\leftrightarrow f}{S}\in F[S]$, or equivalently, $(\frac{x}{S}\vee\frac{e^*}{S})\leftrightarrow\frac{f}{S}\in F[S]$. From $e\in S\cap B(L)$, we have $\frac{e^*}{S}=\frac{o}{S}$ by Lemma 4.4. Consequently, \[\frac{x}{S}\leftrightarrow\frac{f}{S}=(\frac{x}{S}\vee\frac{0}{S})\leftrightarrow\frac{f}{S}=(\frac{x}{S}\vee\frac{e^*}{S})\leftrightarrow\frac{f}{S}\in F[S].\]
		From $f\in B(L)$, we have $\frac{f}{S}\in  B(L[S])$, and so $F[S]$  has BLP in $L[S]$ by Theorem 3.3.
	\end{proof}

	\begin{teorema}\label{b6}
		Let $F$ be a proper filter and $S$ be a $\wedge-$closed system of  a residuated lattice $L$. Then we have a residuated lattice isomorphism \[ L/(F\vee[S\cap  B(L)))\cong L[S]/F[S].\]
	\end{teorema}
	\begin{proof}
		Let  $\psi: L\to L[S]/F[S]$ be defined by $\psi(x):=\frac{x}{S}/F[S]$. Clearly $\psi $ is a surjective morphism of residuated lattices. If $x\in Ker(\psi)$, then $\frac{x}{S}/F[S]=\frac{1}{S}/F[S]$, or equivalently, $\frac{x}{S}\leftrightarrow\frac{1}{S}\in F[S]$. It follows that $\frac{x}{S}=\frac{x\leftrightarrow 1}{S}\in F[S]$. Using Lemma 4.4, there exist $e\in S\cap  B(L)$ and $a\in F$ such that $x\wedge e=a\wedge e$.  Now we have  
		\begin{align*}
			a\odot e=a\wedge e\leq x&\Rightarrow a\leq e\to x\\&\Rightarrow  e\to x\in F.
		\end{align*}
		Now since $e\odot(e\to x)\leq x$ and $e\in S\cap B(L)$, we have $x\in F\vee[S\cap  B(L))$ by Proposition 1.6.
		
		If $y\in F\vee[S\cap  B(L))$, then by Proposition 1.6 there are $a\in F$ and $e\in S\cap B(L)$ such that $a\odot e\leq y$.  We conclude from Proposition 4.2 that $\frac{e}{S}=\frac{1}{S}$, hence we have 
		\begin{align*}
			a\odot e\leq y&\Rightarrow a\leq e\to y\\&\Rightarrow  e\to y\in F\\&\Rightarrow \frac{e\to y}{S}\in F[S]\\&\Rightarrow \frac{e}{S}\to\frac{y}{S}\in F[S]\\&\Rightarrow\frac{1}{S}\to\frac{y}{S}\in F[S]\\&\Rightarrow \frac{y}{S}\in F[S] \\&\Rightarrow \psi(y)=\frac{y}{S}/F[S]=\frac{1}{S}/F[S]\\&\Rightarrow y\in Ker(\psi).
		\end{align*}
		
		Thus $Ker(\psi)=F\vee[S\cap  B(L))$. Therefore, we have $L/(F\vee[S\cap  B(L)))\cong L[S]/F[S]$ by the first isomorphism theorem.
	\end{proof}	
	\begin{lema}\label{a4}
		If $x\in L$ and $e\in B(L)$ such that $e\leq x\vee x^*$, then $e\wedge x, e\wedge x^*, e^*\vee x^{*}, e^*\vee x^{**}\in  B(L) $.
	\end{lema}
	\begin{proof} We shall only prove that $e\wedge x\in  B(L) $. The other cases can similarly be proved.
		Since $e\leq x\vee x^*$, we have 
		\begin{align*}
			e&=e\wedge (x\vee x^*)
			=(e\wedge x)\vee(e\wedge x^*) \\
			&=(e\wedge e\wedge x)\vee(e\wedge x^*)=e\wedge ((e\wedge x)\vee x^*).
		\end{align*}
		Hence $e\leq ((e\wedge x)\vee x^*)$. It follows that
		\begin{align*}
			(e\wedge x)\vee (e\wedge x)^*&=(e\wedge x)\vee (e^*\vee x^*)= ((e\wedge x)\vee  x^*)\vee e^*\geq e\vee e^*=1.
		\end{align*}
		Therefore, $(e\wedge x)\vee (e\wedge x)^*=1$, and so 	$e\wedge x\in  B(L) $ by Proposition 1.2.
	\end{proof}	
	\begin{teorema}\label{b11} 
		Let $F$ be a proper filter of a residuated lattice $L$ such that for each $e\in  B(L)$ we have either $e\in F$ or $e^*\in F$. Set $S_F:=F\cap  B(L)=\{e\mid e\in F\cap  B(L) \}$. Then $S_F$ is a $\wedge-$closed system of $L$ and the residuated lattice $L[S_F]$ is directly indecomposable.
	\end{teorema}
	\begin{proof}
		Since $F$ is a filter, $S_F$ is clearly a $\wedge-$closed system of $L$. If $\frac{x}{S_F}\in B(L[S_F])$, then by Proposition 4.2, there exists $e\in S_F$ such that $e\leq x\vee x^*.$ By Proposition 4.2 and Lemma 4.10, we have $\frac{e}{S_F}=\frac{1}{S_F}$ and $x\wedge e\in B(L)$. It follows that \[\frac{x}{S_F}=\frac{x}{S_F}\wedge\frac{1}{S_F}=\frac{x}{S_F}\wedge\frac{e}{S_F}=\frac{x\wedge e}{S_F}.\]
		By assumption, we have either $x\wedge e\in S_F$ or $(x\wedge e)^*\in S_F$. If $x\wedge e\in S_F$, then $\frac{x}{S_F}=\frac{e\wedge x}{S_F}=\frac{1}{S_F}$, and if $(x\wedge e)^*\in S_F$, then $\frac{x}{S_F}=\frac{e\wedge x}{S_F}=\frac{0}{S_F}$ by Proposition 4.2 and Lemma 4.4. Therefore, $ B(L[S_F])=\{\frac{0}{S_F}, \frac{1}{S_F}\}$, or equivalently, the residuated lattice $L[S_F]$ is directly indecomposable.
	\end{proof}
	
	For a filter $F$ of a residuated lattice $L$, set $F':=[F\cap B(L))=[\{e\mid e\in F\cap B(L) \})$. If $F=L$, then $0\in F\cap B(L)$ and so $F'=L$. Thus $F'=L$ if and only if $F=L$.
	\begin{lema}\label{b13}
		Let $F, G\in    Filt(L)$ for a residuated lattice $L$. Then  $F'\vee G'=L$ if and only if there exists $e\in B(L)$ such that $e\in F$ and $e^*\in G$.
	\end{lema}
	\begin{proof}
		Assume that $F'\vee G'=L$. Thus $0\in F'\vee G'=L$, and so by Proposition 1.6
		there are $a\in F'$ and $b\in G'$ such that $a\odot b=0$. Since $a\in F'$ and $b\in G'$, there are $e\in F\cap  B(L)$ and $f\in G\cap B(L)$ such that $e\leq a$ and $f\leq b$ by Proposition 1.6. Hence $e\odot f\le a\odot b=0$, and so $e\odot f=0$. Thus $f\leq e^*$, and since $f\in G$, we have $e^*\in G$. Therefore, $e\in F$ and $e^*\in G$. The converse is clear.
	\end{proof}

	The following theorem provides some additional conditions that are equivalent to a filter having  BLP.
	\begin{teorema}\label{b14}
		Let $F$ be a proper filter of a residuated lattice $L$. Then the following statements are equivalent:
		
		\begin{enumerate}
			\item $F$ has BLP;
			\item For each maximal (or prime) filter $M$ of $L$, the filter $F[S_M]$   has BLP in $L[S_M]$; 
			
			\item For each maximal (or prime) filter $M$ of $L$, if $F[S_M]\not=L[S_M]$, then $F[S_M]$  is a   pseudo-irreducible filter  of $L[S_M]$; 
			
			\item  For each maximal (or prime) filter $M$ of $L$, if $F\cap (L\setminus S_M)=\varnothing$, then  $F\vee [S_M)$  is a   pseudo-irreducible filter  of $L$;
			
			\item For each maximal (or prime) filter $M$ of $L$, if $F\vee [ S_M)\not=L$ and $e^*\vee x\vee x^*\in F$ for some $e\in S_M$ and $x\in L$, then there exists $f\in S_M$ such that we have either $f^*\vee x\in F$ or $f^*\vee x^*\in F$.
		\end{enumerate}
	\end{teorema}
	
	\begin{proof} 
		$(1)\Rightarrow (2)$. It follows by Proposition 4.8.\\
		$(2)\Rightarrow (3)$. Let $F[S_M]$ be a proper filter of the residuated lattice $L[S_M]$.  By assumption, $F[S_M]$  has BLP. Now by Theorem 4.11, the residuated lattice $L[S_M]$ is directly indecomposable, and hence   $F[S_M]$  is a   pseudo-irreducible filter  of $L[S_M]$ by Proposition 3.6.  
		
		$(3)\Rightarrow (1)$. Let $x\vee x^*\in F$. Set $T:=([x)\to F)'\vee([x^*)\to F)'$. Let $M\in   MaxF(L) $ be such that $T\subseteq M$. If $F[S_M]=L[S_M]$, then by Proposition 4.5 there exists $e\in S_M\subseteq M$ such that $e^*\in F$. It follows that $e^*\in T\subseteq M$ since $F'\subseteq T$, which is a contradiction. Hence, $F[S_M]\not=L[S_M]$, and so $F[S_M]$ is a   pseudo-irreducible filter  of $L[S_M]$ by our assumption. From $x\vee x^*\in F$, we have $\frac{x}{S_M}\vee\frac{x^*}{S_M}=\frac{x\vee x^*}{S_M}\in F[S_M]$. Using Theorem 2.8, we have either $\frac{x}{S_M}\in F[S_M]$ or $\frac{ x^*}{S_M}\in F[S_M]$. If  $\frac{x}{S_M}\in F[S_M]$, then by Lemma 4.4 there exist $a\in F$ and $e\in S_M$ such that $x\wedge e=a\wedge e$. Thus we have 
		\begin{align*}
			a\wedge e\leq x\Rightarrow & a\odot e\leq x\\\Rightarrow &a\leq e\to x\\\Rightarrow& a\leq e^*\vee x.
		\end{align*}
		Now since $a\in F$, we have $e^*\vee x\in F$. Hence by Proposition 1.6, we have 
		\begin{align*}
			e^*\vee x\in F&\Rightarrow [e^*\vee x)\subseteq F\Rightarrow[x)\cap [e^*) \subseteq F\\&\Rightarrow e^*\in [x)\to F \\&\Rightarrow e^*\in ([x)\to F)'.
		\end{align*}
		But from $([x)\to F)'\subseteq T\subseteq M$, we have $e^*\in M$, which is a contradiction. Similarly, if $\frac{x^*}{S_M}\in F[S_M]$, we have a contradiction. Therefore, for each maximal filter $M$ of $L$, we have $T\not\subseteq M$, or equivalently, $T=L$. Thus by Lemma 4.12, there exists $e\in B(L)$ such that $e\in([x)\to F)'$ and $e^*\in ([x^*)\to F)'$. Consequently, $e\in([x)\to F)$ and $e^*\in ([x^*)\to F)$, that is,  $[x\vee e)=[x)\cap[e)\subseteq F$ and $[x^*\vee e^*)=[x^*)\cap[e^*)\subseteq F$. Hence, $x\vee e, x^*\vee e^*\in F$, or equivalently, $x\leftrightarrow e^*\in F$. Therefore, $F$ has BLP by Theorem 3.3.
		
		$(3)\Leftrightarrow (4)$. First of all, note that  for each $e\in  B(L)$ and $M\in  SpecF(L)$ we have either ($e\in M$ and $e^*\not\in M$)  or ($e^*\in M$ and $e\not\in M$). Thus, by Proposition 4.5 $F[S_M]\not=L[S_M]$ if and only if $F\cap (L\setminus S_M)=\varnothing$. Now since $S_M\subseteq  B(L)$, we have  $L[S_M]/F[S_M]\cong L/(F\vee[S_M))$ by Theorem 4.9. Using Theorem 2.8, we have $F[S_M]$  is a   pseudo-irreducible filter  of $L[S_M]$ if and only if the residuated lattice $L[S_M]/F[S_M]$ is non-trivial and directly indecomposable  if and only if the residuated lattice $L/(F\vee[S_M))$ is non-trivial and directly indecomposable if and only if $F\vee[S_M)$  is a   pseudo-irreducible filter  of $L$.
		
		$(3)\Leftrightarrow (5)$. It follows from Proposition 4.7 and the fact that $F[S_M]\not=L[S_M]$ if and only if $F\vee [ S_M)\not=L$.
	\end{proof}

	\begin{propozitie}\label{b15}
		Let $F$ be a filter of a residuated lattice $L$ and $X\subseteq  B(L)$. Then if $F$ has BLP, then $F\vee[X)$ has BLP.
	\end{propozitie}
	\begin{proof}
		Let $M\in   MaxF (L)$. If $X\not\subseteq M$, then there is $e\in X\setminus M$. So $e\in L\setminus S_M$, and thus $e\in(F\vee[X))\cap(L\setminus S_M)$. Hence, $(F\vee[X))\cap(L\setminus S_M)\not=\varnothing$ if $X\not\subseteq M$.
		
		Now assume that $X\subseteq M$. Then $X\subseteq S_M$. Hence $(F\vee[X))\vee[S_M)=F\vee[S_M)$. Now since $F$ has BLP, $(F\vee[X))\vee[S_M)=F\vee[S_M)$ is a pseudo-irreducible filter of $L$ when $(F\vee[X))\cap(L\setminus S_M)=\varnothing$	 
		by Theorem 4.13. Consequently,  $F\vee[X)$ has BLP by by Theorem 4.13.
	\end{proof}

	The following result is a useful consequence of Proposition 4.14 by setting $F=\{1\}$.
	\begin{corolar}\label{b16}
		Let  $X\subseteq  B(L)$. Then the filter $[X)$ has BLP.
	\end{corolar}

	
	\section{Boolean Lifting Property and the radical of filters}\label{4}
	Recall that for a proper filter  $F$  of a residuated lattice  $L$, the intersection of all maximal filters of  $L$  containing  $F$  is called \textit{the radical of  $F$} and is denoted by  $ Rad(F)$, that is,  $ Rad(F) = \bigcap V_{Max}(F)$. Clearly,  $ Rad(F)$  is a filter of $L$,  $F \subseteq  Rad(F) $, and $ Rad( Rad(F))= Rad(F)$. In the special case,  $ Rad(\{1\})$  is the intersection of all maximal filters of  $L $, which is always denoted by  $ Rad(L) $ and is called \textit{the radical of $L$.}

	In this section, we want to consider the relation between a proper filter and its radical from the viewpoint of  Boolean lifting property.  We start with the following lemmas.
	
	\begin{lema}\label{b17}
		Let $L$ be a residuated lattice. Then for $F, G, H\in   Filt(L)$ if $F=G\cap H$ and $G\vee H=L$, then $ Rad(F)= Rad(G)\cap Rad(H)$ and $ Rad(G)\vee Rad(H)=L$.
	\end{lema}
	\begin{proof}
		Since $F=G\cap H$ and $G\vee H=L$, we have $V(F)=V(G)\cup V(H)$, $V(G)\cap V(H)=\varnothing$,  $V_{Max}(F)=V_{Max}(G)\cup V_{Max}(H)$, and $V_{Max}(G)\cap V_{Max}(H)=\varnothing$ by Proposition 1.8. Therefore, $ Rad(F)=\bigcap V_{Max}(F)=(\bigcap V_{Max}(G))\cap(\bigcap V_{Max}(H))= Rad(G)\cap Rad(H)$ and since $G\subseteq  Rad(G)$ and $H\subseteq  Rad(H)$, we have  $ Rad(G)\vee Rad(H)=L$. 
	\end{proof}
	\begin{lema}\label{b18}
		For a proper filter $F$  of a residuated lattice $L$ and $e\in B(L)$, if $e\in Rad(F)$, then $e\in F$.
	\end{lema}
	\begin{proof}
		Let $P\in V(F)$. Thus there is $M\in V_{Max}(F)$ with $P\subseteq M$. 
		From $e\in Rad(F)$, we have $e\in M$. Now since $e\wedge e^*=e\odot e^*=0\in P$, we have either $e\in P$ or $e^*\in P$. If $e^*\in P$, then $e^*\in M$, which is a contradiction. So we deduce that $e\in P$ for each $P\in V(F)$. Therefore, $e\in\bigcap V(F)$. By Proposition 1.7, $F=\bigcap V(F)$. Therefore,  $e\in F$.
	\end{proof}
	\begin{teorema}\label{b19}
		For a filter $F$  of a residuated lattice $L$, if $ Rad(F)$ has BLP, then $F$ has BLP.
	\end{teorema}
	\begin{proof}
		Let $F=G\cap H$ and $G\vee H=L$ for some $G, H\in   Filt(L)$. By Lemma 5.1, $ Rad(F)= Rad(G)\cap Rad(H)$ and $ Rad(G)\vee Rad(H)=L$.  By Theorem 3.3 since  $ Rad(F)$ has BLP, there exists $e\in B(L)$ such that $e\in Rad(G)$ and $e^*\in Rad(H)$.   Lemma 5.2 now implies that $e\in\ G$ and $e^*\in H$. Consequently, $F$ has BLP by Theorem 3.3. 
	\end{proof} 

	In [10, Open question 3.4], the authors posed  the question \lq\lq Can sufficient, or even necessary and sufficient conditions be provided for a residuated
	lattice $L$ to be such that $ Rad(L)$ has BLP?". In the following, we  answer this question.

	\begin{definitie}\label{tprd}
		We say that a residuated lattice $L$ has \textit{transitional property of radicals decomposition (TPRD)}, whenever for $F, G, H\in   Filt(L)$ if $ Rad(F)=G\cap H$ and $G\vee H=L$, then there exist  $G_0, H_0\in   Filt(L)$ such that $G_0\subseteq G$, $H_0\subseteq H$, $ Rad(G_0)=G$, $ Rad(H_0)=H$, $F=G_0\cap H_0$, and $G_0\vee H_0=L$.
	\end{definitie}	
	\begin{teorema}\label{b22}
		Let $L$ be a  residuated lattice that has TPRD and $F\in   Filt(L)$. Then $F$ has BLP if and only if $ Rad(F)$ has BLP. In particular, $ Rad(L)$ has BLP. 
	\end{teorema}
	\begin{proof}
		Assume that $F$ has BLP and   $ Rad(F)=G\cap H$ and $G\vee H=L$ for some $ G, H\in   Filt(L)$. By assumption, there exist  $G_0, H_0\in   Filt(L)$ such that $G_0\subseteq G$, $H_0\subseteq H$, $ Rad(G_0)=G$, $ Rad(H_0)=H$, $F=G_0\cap H_0$, and $G_0\vee H_0=L$. Now since $F$ has BLP, by Theorem 3.3 there exists $e\in B(L)$ such that $e\in G_0$ and $e^*\in H_0$. Therefore, $e\in G$ and $e^*\in H$, and so $ Rad(F )$ has BLP by Theorem 3.3. The converse follows from Theorem 5.3.
		
		By Example 3.2, the trivial filter $\{1\}$  always has BLP, hence $ Rad(L)= Rad(\{1\})$ has BLP if $L$  has TPRD.
	\end{proof}

	\begin{definitie}
		A residuated lattice $L$ is called \textit{weak MTL-algebra}, whenever for each $x\in L$ we have $(x^*\to x^{**})\vee(x^{**}\to x^*)=1$.
	\end{definitie}
	
	\begin{exemplu}Let $L = \{0, a, b,  c, d, 1\}$ be a lattice whose  Hasse diagram is given in the following figure:
		
		\begin{center}
			$
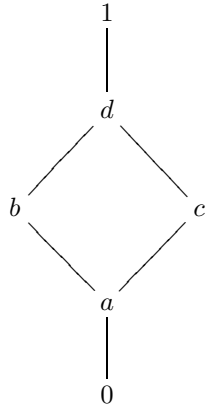
\xymatrix{
				& & 1\ar@{-}[d]\\&& d\ar@{-}[ld]\ar@{-}[rd]\\
				& b \ar@{-}[rd] & & c\ar@{-}[ld]  & \\&& a
				&\\
				&& 0\ar@{-}[u]}$
		\end{center}
		
		Consider the operations  $\odot$ and  $\rightarrow$ given by the
		following tables (see [17, Example 1, pp. 240]):
		
		\begin{center}	\begin{tabular}{r|l l l l l l  }
				$\odot$ & 0  & $a$ & $b$   & $c$ & $d$ & $1$   \\\hline
				0 & $0$ & $0$ &$0$  & $0$ & $0$   & $0$ \\
				
				$a$ & $0$ & $0$ &$0$  & $0$ & $0$   & $a$\\
				$b$ & $0$ & $0$ &$b$  & $0$ & $b$  & $b$ \\
				$c$ & $0$ & $0$ &$0$  & $c$ & $c$  & $c$\\
				$d$ & $0$ & $0$ &$b$  & $c$ & $d$   & $d$ \\
				1 & 0  & $a$ &$b$  & $c$ & $d$ & $1$ 
			\end{tabular}
			\vspace*{0.5cm}
			\begin{tabular}{r|l l l l l l  }
				$\rightarrow$  &0 & $a$ & $b$   & $c$ & $d$ & $1$   \\\hline
				0 & $1$ & $1$ &$1$  & $1$ & $1$   & $1$ \\
				
				$a$ & $d$ & 1 &$1$  & $1$ & $1$  & $1$\\
				$b$ & $c$ & $c$ &$1$  & $c$ & $1$   & $1$ \\
				$c$ & $b$ & $b$ &$b$  & $1$ & $1$   & $1$\\
				$d$ & $a$ & $a$ &$b$  & $c$ & $1$   & $1$ \\
				1 & 0  & $a$ &$b$  & $c$ & $d$ & $1$ 
			\end{tabular}

			\vspace{0.5cm}
		\end{center}
		Since $(b^*\to b^{**})\vee(b^{**}\to b^*)=(c\to b)\vee(b\to c)=b\vee c=d\not=1$,  $L$ is not a weak MTL-algebra. Also, note that an easy computation shows that $L$ is an involution. Hence, an involution residuated lattice need not be  a weak MTL-algebra.
	\end{exemplu} 
	\begin{propozitie}\label{semig}[2, Proposition 9]
		The following conditions are
		equivalent for a residuated lattice $L$.
		\begin{enumerate}
			\item $L$ is a semi-$G$-algebra;
			\item For every $x \in L$,  $x \wedge x^*= 0$.
		\end{enumerate}
	\end{propozitie}
	
	\begin{propozitie}\label{yyyy}
		Every semi-$G$-algebra that is a De Morgan residuated lattice (or equivalently,  Stonean residuated lattice, see [2, Proposition 8 and Corollary 1]) is a weak MTL-algebra.
	\end{propozitie}
	\begin{proof}
		Let $x\in L$. Then by Propositions 1.2 and 5.8, we have 
		\begin{align*}
			(x^*\to x^{**})\vee(x^{**}\to x^*)=&(x^*\to (x^{*}\to 0))\vee(x^{**}\to (x\to 0))\\=&((x^*)^2\to 0)\vee((x^{**}\odot x)\to 0)\\=&((x^*)^2)^*\vee((x^{**}\odot x)^*)\\\geq & x^{**}\vee x^*=(x^*\wedge x)^*=0^*=1.
		\end{align*}
		Therefore, $(x^*\to x^{**})\vee(x^{**}\to x^*)=1$ and so $L$ is a weak MTL-algebra.
	\end{proof}

	Other examples of weak MTL-algebras are Boolean algebras, MV-algebras,
	BL-algebras, MTL-algebras. The class of weak MTL-algebras strictly contains the class of  MTL-algebras and Stonean residuated lattices.	
	
	In the following examples and Example 5.15, we show that all assumptions of Proposition 5.9 are necessary. Also, the  converse of Proposition 5.9 need not be true.
	\begin{exemplu}Let $L = \{0, n, a, b, , c, d, 1\}$ be a lattice whose  Hasse diagram is given in the following figure:
		
		\begin{center}
			$
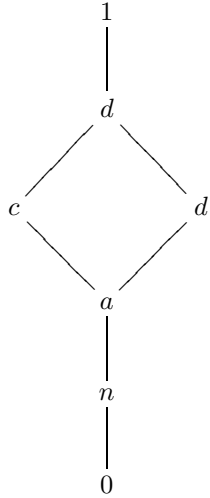
\xymatrix{
				& & 1\ar@{-}[d]\\&& d\ar@{-}[ld]\ar@{-}[rd]\\
				& c \ar@{-}[rd] & & d\ar@{-}[ld]  & \\&& a
				&\\
				&  &n\ar@{-}[u] &&\\
				&& 0\ar@{-}[u]}$
		\end{center}
		
		Consider the operations  $\odot$ and  $\rightarrow$ given by the
		following tables (see [12, Example 3.2]):
		\begin{center}
			\begin{tabular}{r|l l l l l l l }
				$\odot$ & 0 & $n$ & $a$ &$b$  & $c$ & $d$ & $1$   \\\hline
				0 & $0$ & $0$ &$0$  & $0$ & $0$ & $0$  & $0$ \\
				$n$& $0$ & $0$ &$0$  & $n$ & $0$ & $n$  & $n$\\
				$a$ & $0$ & $0$ &$a$  & $a$ & $a$ & $a$  & $a$\\
				$b$ & $0$ & $n$ &$a$  & $b$ & $a$ & $b$  & $b$ \\
				$c$ & $0$ & $0$ &$a$  & $a$ & $c$ & $c$  & $c$\\
				$d$ & $0$ & $n$ &$a$  & $b$ & $c$ & $d$  & $d$ \\
				1 & 0 &$n$ & $a$ &$b$  & $c$ & $d$ & $1$ 
			\end{tabular}
			\vspace*{0.5cm}
			\begin{tabular}{r|l l l l l l l }
				$\rightarrow$  &0& $n$ & $a$ &$b$  & $c$ & $d$ & $1$   \\\hline
				0 & $1$ & $1$ &$1$  & $1$ & $1$ & $1$  & $1$ \\
				$n$& $c$ & $1$ &$1$  & $1$ & $1$ & $1$  & $1$\\
				$a$ & $n$ & $n$ &$1$  & $1$ & $1$ & $1$  & $1$\\
				$b$ & $0$ & $n$ &$c$  & $1$ & $c$ & $1$  & $1$ \\
				$c$ & $n$ & $n$ &$b$  & $b$ & $1$ & $1$  & $1$\\
				$d$ & $0$ & $n$ &$a$  & $b$ & $c$ & $1$  & $1$ \\
				1 & 0 &$n$ & $a$ &$b$  & $c$ & $d$ & $1$ 
			\end{tabular}

			\vspace{0.5cm}
		\end{center}
		An easy computation shows that  $(x^*\to x^{**})\vee(x^{**}\to x^*)=1$ for each $x\in L$, that is,   $L$ is a weak MTL-algebra. But since $(b\to c)\vee(c\to b)=c\vee b=d\not=1$,  $L$ is not an MTL-algebra. Also, since $(n^2)^*=0^*=1\not=c=n^*$, $L$ is neither semi-$G$-algebra nor Stonean residuated lattice, see [2, Proposition 8 and Corollary 1] for more details. Also, since $a^{**}=c\not=a$, $L$ is not an involution residuated lattice.
	\end{exemplu}
	
	\begin{exemplu}Let $L = \{0, n, a, b, c, d, e, f, m, 1\}$ be a lattice whose  Hasse diagram is given in the following figure:
		
		\begin{center}
			$
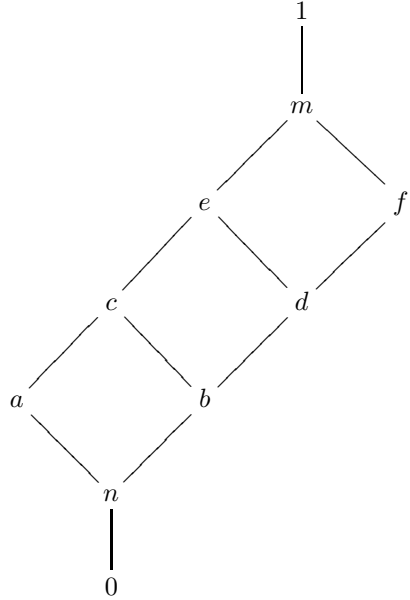
\xymatrix{
				& && & 1\ar@{-}[d]\\&&&& m\ar@{-}[rd]\\
				& && e\ar@{-}[ru]& & f\  & \\&&c\ar@{-}[ru]& &d\ar@{-}[ru]\ar@{-}[lu] 
				&\\&a\ar@{-}[ru]&&b\ar@{-}[ru]\ar@{-}[lu]& 
				&\\
				&&n\ar@{-}[ru]\ar@{-}[lu]&& \\&&0\ar@{-}[u]&&}$
		\end{center}
		
		Consider the operations  $\odot$ and  $\rightarrow$ given by the
		following tables (see [15, Example 4.2]):
		
		\begin{center}
			
			\begin{tabular}{r|l l l l l l l l l l }
				$\odot$ & 0 & $n$ & $a$ & $b$   & $c$ & $d$ & $e$& $f$& $m$ & $1$   \\\hline
				0	& 0 & $0$ & $0$ & $0$   & $0$ & $0$ & $0$& $0$& $0$ & $0$ \\
				
				$n$& 0 & $0$ & $0$ & $0$   & $0$ & $0$ & $0$& $0$& $0$ & $n$\\
				$a$	& 0 & $0$ & $a$ & $0$   & $a$ & $0$ & $a$& $0$& $a$ & $a$ \\
				$b$	& 0 & $0$ & $0$ & $0$   & $0$ & $0$ & $0$& $b$& $b$ & $b$\\
				$c$	& 0 & $0$ & $a$ & $0$   & $a$ & $0$ & $a$& $b$& $c$ & $c$ \\
				$d$	& 0 & $0$ & $0$ & $0$   & $0$ & $b$ & $b$& $d$& $d$ & $d$\\
				$e$	& 0 & $0$ & $a$ & $0$   & $a$ & $b$ & $c$& $d$& $e$ & $e$\\
				$f$	& 0 & $0$ & $0$ & $b$   & $b$ & $d$ & $d$& $f$& $f$ & $f$\\
				$m$	& 0 & $0$ & $a$ & $b$   & $c$ & $d$ & $e$& $f$& $m$ & $m$\\
				1	& 0 & $n$ & $a$ & $b$   & $c$ & $d$ & $e$& $f$& $m$ & $1$ 
			\end{tabular}
			\vspace{0.5cm}
			\begin{tabular}{r|l l l l l l l l l l }
				$\rightarrow$  & 0 & $n$ & $a$ & $b$   & $c$ & $d$ & $e$& $f$& $m$ & $1$   \\\hline
				0	& 1 & $1$ & $1$ & $1$   & $1$ & $1$ & $1$& $1$& $1$ & $1$ \\	
				$n$	& $m$ & $1$ & $1$ & $1$   & $1$ & $1$ & $1$& $1$& $1$ & $1$\\
				$a$	& $f$ & $f$ & $1$ & $f$   & $1$ & $f$ & $1$& $f$& $1$ & $1$ \\
				$b$	& $e$ & $e$ & $e$ & $1$   & $1$ & $1$ & $1$& $1$& $1$ & $1$\\
				$c$	& $d$ & $d$ & $e$ & $f$   & $1$ & $f$ & $1$& $f$& $1$ & $1$ \\
				$d$	& $c$ & $c$ & $c$ & $e$   & $e$ & $1$ & $1$& $1$& $1$ & $1$\\
				$e$	& $b$ & $b$ & $c$ & $d$   & $e$ & $f$ & $1$& $f$& $1$ & $1$\\
				$f$	& $a$ & $a$ & $a$ & $c$   & $c$ & $e$ & $e$& $1$& $1$ & $1$\\
				$m$	& $n$ & $n$ & $a$ & $b$   & $c$ & $d$ & $e$& $f$& $1$ & $1$\\
				1	& 0 & $n$ & $a$ & $b$   & $c$ & $d$ & $e$& $f$& $m$ & $1$
			\end{tabular}
			\vspace{0.5cm}
		\end{center}
		By [15, Example 4.2] and an easy argument we deduce that  $L$ is a De Morgan residuated lattice, but $L$ is not a weak MTL-algebra since $(a^*\to a^{**})\vee(a^{**}\to a^*)=(f\to a)\vee(a\to f)=a\vee f=m\not=1$. Also, $L$ is not a semi-$G$-algebra since $(b^2)^*=0^*=1\not=e=b^*$. 
	\end{exemplu}

	\begin{propozitie}\label{3.18} For a residuated lattice $L$ let
		$F_{1},..., F_{n}$
		be  $n$ proper pairwise comaximal filters of $L$ (that is, $F_i\vee F_j=L$ for each $i\not=j$) and $G_{1},..., G_{n}$ be $n $ filters of $L$  
		such that  $F_i\subseteq G_i$ for each $1\leq i\leq n$. If
		$\cap_{i=1}^nF_i=\cap_{i=1}^nG_i$, then $F_i=G_i$ for all for $1\leq i\leq n$.
		
		\begin{proof} Since $(   Filt(L), \subseteq)$ is a distributive lattice, we have  $F_1\vee (\cap_{i=2}^nF_i)=\cap_{i=2}^n(F_1\vee F_i)=\cap_{i=2}^nL=L$.
			Then by Proposition 1.2, there exist $a\in F_1$ and $ b\in \cap_{i=2}^nF_i$ such that $a\odot b=0$.
			Let $x\in G_1$ be arbitrary. Then $x\vee b\in G_1\cap(\cap_{i=2}^nF_i)
			\subseteq \cap_{i=1}^nG_i=\cap_{i=1}^nF_i\subseteq	F_1$. Also, $x\vee a\in F_1$. Using Proposition 1.2, we have 
			$x=x\vee 0=x\vee (a\odot b)\geq (x\vee a)\odot(x\wedge b)\in F_1$. Hence $x\in F_1$. Consequently, $F_1=G_1$.
			Similarly, we can prove $F_i=G_i$ for $i=2,..., n$.
		\end{proof}
	\end{propozitie}
	
	\begin{corolar}\label{edn}
		Let $F, G, H\in   Filt(L)$. If $ Rad(F)=G\cap H$ and $G\vee H=L$, then $G= Rad(G)$ and $H= Rad(H)$, or equivalently, $G=\bigcap V_{Max}(G)$ and $H=\bigcap V_{Max}(H)$.
	\end{corolar}
	\begin{proof} From Lemma 5.1 and the fact that $ Rad( Rad(F))= Rad(F)$,
		we have $G\cap H= Rad(G)\cap  Rad(H)$. Since $G\vee H=L$,  by Proposition 5.12 and the fact that $G\subseteq  Rad(G)$ and $H\subseteq Rad(H)$, we have $G= Rad(G)$ and $H= Rad(H)$, or equivalently, $G=\bigcap V_{Max}(G)$ and $H=\bigcap V_{Max}(H)$.
	\end{proof}
	
	In [8, Corollary 6.16], it was shown that $ Rad(L)$ has BLP for MV-algebras and BL-algebras; in the following proposition, we extend this result.
	\begin{teorema}\label{b20} 
		Every weak MTL-algebra (e.g., Boolean-algebra, MV-algebra, BL-algebra, MTL-algebra, Stonean residuated lattice) $L$ has TPRD. In particular,  for every weak MTL-algebra (e.g., Boolean-algebra, MV-algebra, BL-algebra, MTL-algebra, Stonean residuated lattice) $L$, $ Rad(L)$ has BLP.
	\end{teorema}
	\begin{proof} Let $F, G, H\in   Filt(L)$ such that $ Rad(F)=G\cap H$ and $G\vee H=L$.
		Since $G\vee H=L$, there exist $x\in G$ and $y\in H$ such that $x\odot y=0$ by Proposition 1.6. Thus $y\leq x^*$, and so $x^*\in H$. Also, since $x\leq x^{**}$, we have $x^{**}\in G$. 
		
		Let $M\in V_{Max}(F)$. From $M\in V_{Max}(F)$ we conclude that $G\cap H= Rad(F)\subseteq M$. Hence,  we have either $G\subseteq M$ or $H\subseteq M$ from the fact that every maximal filter is prime. By assumption we have $(x^*\to x^{**})\vee(x^{**}\to x^*)=1$, and so $(x^*\to x^{**})\vee(x^{**}\to x^*)\in M$. It follows that either $x^*\to x^{**}\in M$ or $x^{**}\to x^{*}\in M$.
		
		Now if  $G\subseteq M$ and  $x^{**}\to x^{*}\in M$, then from $x^{**}\in G$ we have $x^*, x^{**}\in M$, which is impossible. Also, if $H\subseteq M$ and  $x^{*}\to x^{**}\in M$, then from $x^{*}\in H$ we have $x^*, x^{**}\in M$, which is impossible. Therefore, $G\subseteq M$ if and only if  $x^*\to x^{**}\in M$, and $H\subseteq M$ if and only if $x^{**}\to x^{*}\in M$, that is, $V_{Max}(G)=V_{Max}(x^*\to x^{**})\cap V_{Max}(F)$ and $V_{Max}(H)=V_{Max}(x^{**}\to x^{*})\cap  V_{Max}(F)$.
		
		Set $A:=\{P\in V(F)\mid x^*\to x^{**}\in P\}$ and $B:=\{P\in V(F)\mid x^{**}\to x^{*}\in P\}$. Let $P\in V(F)$. From  $(x^*\to x^{**})\vee(x^{**}\to x^*)=1\in P$, we have either $x^*\to x^{**}\in P$ or $x^{**}\to x^*\in P$.  Consequently, we have $V(F)=A\cup B$.
		If $P\in A\cap B$, then $x^*\to x^{**}, x^{**}\to x^{*} \in P$. Since $P$ is a proper filter of $L$ containing $F$, there exists $M\in V_{Max}(F)$ such that $P\subseteq M$, and so $x^*\to x^{**}, x^{**}\to x^{*} \in M$, that is, $M\in V_{Max}(x^*\to x^{**})\cap V_{Max}(x^{**}\to x^{*})\cap V_{Max}(F)=V_{Max}(G)\cap V_{Max}(H)=\varnothing$, which is impossible. Therefore, 
		$A\cap B=\varnothing$.
		
		Set $G_0:=\bigcap A$ and $H_0:=\bigcap B$.			
		Clearly, $G_0\subseteq G$, $H_0\subseteq H$. By the above argument, $V_{Max}(G)=V_{Max}(x^*\to x^{**})\cap V_{Max}(F)=A\cap  MaxF (L)$ and $V_{Max}(H)=V_{Max}(x^{**}\to x^{*})\cap V_{Max}(F)=B\cap  MaxF (L)$. Thus, by Corollary 5.13, $ Rad(G_0)=\bigcap(A\cap  MaxF (L))=\bigcap V_{Max}(G)=G$ and  $ Rad(H_0)=\bigcap(B\cap  MaxF (L))=\bigcap V_{Max}(H)=H$. Also, from  $V(F)=A\cup B$ and $A\cap B=\varnothing$, we have  $F=G_0\cap H_0$, and $G_0\vee H_0=L$. Therefore, $L$ has TPRD.  The rest of the theorem follows from Theorem 5.5.
	\end{proof}

	\begin{exemplu}\label{wgl}Let $L = \{0, a, b,  c, d, 1\}$ be a lattice whose  Hasse diagram is given in the following figure:
		
		\begin{center}
			$
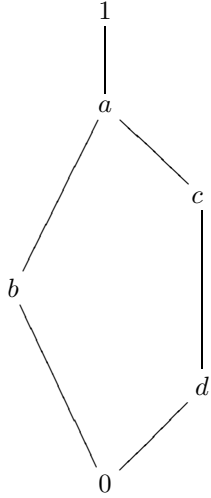
\xymatrix{
				& & 1\ar@{-}[d]\\&& a\ar@{-}[ldd]\ar@{-}[rd]\\
				&  & & c\  & \\& b& 
				&\\& & 
				&d\ar@{-}[uu]\ar@{-}[ld]\\
				&& 0\ar@{-}[luu]}$
		\end{center}
		
		Consider the operations  $\odot$ and  $\rightarrow$ given by the
		following tables (see [8, Example 4.16]):
		
		\begin{center}
			
			\begin{tabular}{r|l l l l l l  }
				$\odot$ & 0  & $a$ & $b$   & $c$ & $d$ & $1$   \\\hline
				0 & $0$ & $0$ &$0$  & $0$ & $0$   & $0$ \\
				
				$a$ & $0$ & $a$ &$b$  & $d$ & $d$   & $a$\\
				$b$ & $0$ & $b$ &$b$  & $0$ & $0$  & $b$ \\
				$c$ & $0$ & $d$ &$0$  & $d$ & $d$  & $c$\\
				$d$ & $0$ & $d$ &$0$  & $d$ & $d$   & $d$ \\
				1 & 0  & $a$ &$b$  & $c$ & $d$ & $1$ 
			\end{tabular}
			\vspace{0.5cm}
			\begin{tabular}{r|l l l l l l  }
				$\rightarrow$  &0 & $a$ & $b$   & $c$ & $d$ & $1$   \\\hline
				0 & $1$ & $1$ &$1$  & $1$ & $1$   & $1$ \\
				
				$a$ & $0$ & 1 &$b$  & $c$ & $c$  & $1$\\
				$b$ & $c$ & $1$ &$1$  & $c$ & $c$   & $1$ \\
				$c$ & $b$ & $1$ &$b$  & $1$ & $a$   & $1$\\
				$d$ & $b$ & $1$ &$b$  & $1$ & $1$   & $1$ \\
				1 & 0  & $a$ &$b$  & $c$ & $d$ & $1$ 
			\end{tabular}
			\vspace{0.5cm}
		\end{center}
		By [8, Example 4.16], $ Rad(L)$ has BLP (actually, every filter of  $L$ has BLP), but $L$ is neither  a chain, nor local, nor a Boolean algebra,  nor a $G$-algebra. Also, since $(b^*\to b^{**})\vee(b^{**}\to b^*)=(c\to b)\vee(b\to c)=b\vee c=a\not= 1$ and $b^*\vee b^{**}=c\vee c^*=c\vee b=a\not=1$, we have $L$ is neither weak MTL-algebra, nor Stonean algebra. An easy computation shows that $L$ is a semi-$G$-algebra.
	\end{exemplu}
	
	We end this section with a topological characterization for residuated lattices whose radical has BLP. This characterization is a topological answer  to [10, Open question 3.4].
	\begin{teorema}\label{b24}
		Let $L$ be a residuated lattice. Then $ Rad(L)$ has BLP if and only if \[   Clop(  MaxF (L))=\{V_{Max}(e)\mid e\in  B(L)\},\] where $   Clop$ is the set of all closed and open subsets of $  MaxF(L)$ with respect to the Stone topology.
	\end{teorema}
	\begin{proof}
		$\Rightarrow).$ Clearly, $\{V_{Max}(e)\mid e\in  B(L)\}\subseteq    Clop(  MaxF (L))$. Let $T$ be a clopen subset of $  MaxF (L)$. Then there exist $A, B\in   Filt(L)$ such that $T=V_{Max}(A)$ and $T^c:=  MaxF(L)\setminus T=V_{Max}(B) $. Set $F:=\bigcap T$ and $G:=\bigcap  T^c$. Hence, $ Rad(L)=F\cap G$. If $M\in  MaxF (L)$ such that $F, G\subseteq M$, then $A\subseteq F\subseteq M$ and $B\subseteq G\subseteq M$. Consequently, $M\in T\cap T^c=\varnothing$, which is impossible. Therefore, $F\vee G=L$. Now by assumption and Theorem 3.3, there exists $e\in B(L)$ such that $e\in F$ and $e^*\in G$. It follows that $T\subseteq V_{Max}(e)$ and $T^c\subseteq V_{Max}(e^*) $. From $  MaxF (L)=T\cup T^c$, $  MaxF (L)=V_{Max}(e)\cup V_{Max}(e^*)$, $T\cap T^c=\varnothing$, and $V_{Max}(e)\cap V_{Max}(e^*)=\varnothing$, we can assert that $T= V_{Max}(e)$ and $T^c= V_{Max}(e^*) $.

		$\Leftarrow).$ Let $x\vee x^*\in  Rad(L)$ for some $x\in L$. Hence for each $M\in  MaxF (L)$, we have $x\vee x^*\in M$. Thus for each $M\in  MaxF (L)$, we have either $x\in M$ or $x^*\in M$, equivalently,
		$  MaxF (L)=V_{Max}(x)\cup V_{Max}(x^*)$. Now since $x\odot x^*=0$, we conclude that 
		$V_{Max}(x)\cap V_{Max}(x^*)=\varnothing$. It follows that $V_{Max}(x)$ is a clopen subset of $  MaxF (L)$. Thus by our assumption, there exists $e\in  B(L)$ such that $V_{Max}(x)=V_{Max}(e)$ and $V_{Max}(x^*)=V_{Max}(e^*)$. Now let $M\in   MaxF(L)$. If $M\in V_{Max}(x)=V_{Max}(e)$, then $x, e\in M$, and so $x\leftrightarrow e=(x^*\vee e)\wedge(x\vee e^*)\in M$. If $M\in V_{Max}(x^*)=V_{Max}(e^*)$, then $x^*, e^*\in M$, and so $x\leftrightarrow e=(x^*\vee e)\wedge(x\vee e^*)\in M$. Therefore, for each $M\in   MaxF (L)$ we have $x\leftrightarrow e\in M$, or equivalently, $x\leftrightarrow e\in Rad(L)$. It follows that $ Rad(L)$ has BLP by Theorem 3.3.
	\end{proof}

\end{document}